%% file: AbSurf_Complete.tex
\documentclass{amsart}
\usepackage[utf8]{inputenc}
\usepackage{amsaddr}
\usepackage{latexsym,amssymb}
\usepackage{amsmath}
\usepackage{tabularx}
\usepackage{hyperref}
\usepackage[matrix,arrow,curve]{xy}

\numberwithin{equation}{section}
\newtheorem{theorem}{Theorem}[section]
\newtheorem{lemma}[theorem]{Lemma}
\newtheorem{proposition}[theorem]{Proposition}
\newtheorem{corollary}[theorem]{Corollary}

\theoremstyle{definition}
\newtheorem{definition}[theorem]{Definition}

\newtheorem{remark}[theorem]{Remark}
\newtheorem{example}[theorem]{Example}

\newcommand{\ZZ}{\mathbb{Z}}

\renewcommand{\AA}{\mathbb{A}}
\newcommand{\VV}{\mathbb{V}}

\newcommand{\PP}{\mathbb{P}}
\newcommand{\FF}{\mathcal{F}}

\newcommand{\Q}{\mathcal{Q}}

\newcommand{\OO}{\mathcal{O}}
\newcommand{\II}{\mathcal{I}}
\newcommand{\RR}{\mathcal{R}}

\newcommand{\tr}{\frac{1}{d}\operatorname{tr}}

\newcommand{\Sym}{\operatorname{Sym}}

\renewcommand{\OO}{\mathcal{O}}
\newcommand{\EE}{\mathcal{E}}
\newcommand{\UU}{\mathcal{U}}

\newcommand{\Spec}{\operatorname{Spec}}
\newcommand{\SSpec}{\operatorname{\underline{Spec}}}

\newcommand{\Proj}{\operatorname{Proj}}

\newcommand{\Hom}{\operatorname{Hom}}

\newcommand{\oGr}{\operatorname{oGr}}
\newcommand{\GL}{\operatorname{GL}}
\newcommand{\PGL}{\operatorname{PGL}}
\newcommand{\LL}{\mathcal{L}}
\newcommand{\MM}{\mathcal{M}}
\newcommand{\HHom}{\mathcal{H}om}

\begin{document}

\title[On the Graded Equations of $(1,3)$-Abelian Surfaces]{On the Graded Equations of $(1,3)$-Abelian Surfaces}

\author{Eduardo Dias}
\address{CMUP, Department of Mathematics, University of Porto,  Portugal}
\email{eduardo.dias@fc.up.pt}
\thanks{This research was supported by FCT (Portugal) under the project 
PTDC/MAT-GEO/2823/2014 and by CMUP (UID/MAT/00144/2019),
which is funded by FCT with national (MCTES) and European structural funds through the programs FEDER under the partnership agreement PT2020.}

\maketitle

\begin{abstract}
Let $S$ be an abelian surface over an algebraically closed field $k$ with characteristic different from $2$ and $3$, and $\LL$ a symmetric ample line bundle defining a polarisation of type $(1,3)$. Then the linear system $|\LL|$ defines a covering map $\varphi\colon S\rightarrow \PP^2$ of degree $6$. Furthermore, if $|\LL|$ is base point free, then $\varphi_*\OO_S=\OO_{\PP^2}\oplus\Omega^1_{\PP^2}\oplus\Omega^1_{\PP^2}\oplus\OO_{\PP^2}(-3)$.

Using this decomposition, in this paper we construct the graded coordinate ring of $(S,\LL,\theta)$, where $\theta\colon G(\LL)\xrightarrow{\sim} H(1,3)$ is a level structure of canonical type. As a corollary we prove that the moduli space of such triples is rational.
\end{abstract}

%\begin{acknowledgements}
% 
%\end{acknowledgements}

\section{Introduction}
\input IntroductionAugust2019.tex

\section{Total Space of a finite rank locally free sheaf}\label{Preliminaries}
\input RelativeAffineSpacePn.tex

\section{Cover Homomorphisms}\label{linearalgebra}
\input CoverHomomorphism.tex

\section{Coordinate Ring and Moduli}\label{CoordinateRing}
\input CoordinateRing.tex

\section{Appendix}
\input Appendix.tex

\bibliographystyle{amsalpha}
\bibliography{Lit}

\end{document}

%% file: IntroductionAugust2019.tex
Let $S$ be an abelian surface over an algebraically closed field $k$ with characteristic different from $2$ and $3$, and let $\LL$ be an ample, symmetric line bundle on $S$ defining a polarisation of type $(1,d)$. In this paper our first result is the description of the coordinate ring of $(S,\LL,\theta)$, i.e. 
$$\RR(S,\LL,\theta)=\bigoplus_{n\in\ZZ}H^0(S,\LL^n),$$
when $d=3$, and $\theta\colon G(\LL)\xrightarrow{\sim} H(1,3)$ is a level structure of canonical type. The methods will be used in a follow up paper for $d=4$.  
 
Although a complex torus of dimension $g$ has a simple description, the quotient of a $g$-dimensional complex vector space by a maximal rank lattice, writing explicit equations describing the coordinate ring of one satisfying the Riemann conditions, i.e. an abelian variety, was always elusive for any dimension greater than $1$.
 
Partial results are known. A $(1,4)$-abelian surface is birational to a singular octic in $\PP^3$, see \cite{BLange4}, the well known $(1,5)$-abelian surface is the zero set of a section of the Horrocks-Mumford bundle, see \cite{Mumford5}, and birational to a surface in $\PP^4$ as described in \cite{15embedding}. In \cite{MSchreyer7} the equations and syzygies of the embedding of $(1,7)$-abelian surfaces in $\PP^6$ are described, and in \cite{GrossPopescu} the authors describe the homogeneous ideal of the embedding of $A\times A\subset \PP^{d-1}\times \PP^{d-1}$, where $A$ is a $(1,d)$-abelian surface for $d\geq 10$. All these constructions fail to capture the irregularity of the coordinate ring of an abelian surface.

From an algebraic point of view, irregularity forces the failure of the \textit{Arithmetically Cohen--Macaulay} condition on the coordinate ring $\RR(S,\LL)$. Because of it, one has to replace the Hilbert Series method used to describe the free resolution of the coordinate ring of regular varieties by a different one. For a pair $(S,\LL)$, where $\LL$ defines a polarisation of type $(1,d)$, the linear system $|\LL|$ induces a rational map $\varphi_\LL\colon S\rightarrow \PP=\PP^{d-1}$. Then, under mild conditions, one can use Beilinson’s theorem \cite{Beilinson1} to describe the structure of the ring as an $\OO_{\PP^{d-1}}$-module. 

For the case of a $(1,3)$-abelian surface, assuming that $(S,\LL)\not\simeq(E_1\times E_2,p_1^*\LL_1\otimes p_2^*\LL_2)$, where $E_i$ are elliptic curves and $\LL_i$ are line bundles on $E_i$, then, by \cite[Lemma $10.1.1$]{BLange}, $|\LL|$ has no fixed components and is base point free. We get that $\varphi$ is a surjective map and so a covering map of degree $\LL^2=2d=6$. Moreover, using Beilinson's theorem, it is proven in \cite[$0.4$]{Casab} that the linear system $|\LL|$ defines a morphism $\varphi\colon S\rightarrow \PP^2$ such that 
$$\varphi_*\OO_S=\OO_{\PP^2}\oplus\Omega^1_{\PP^2}\oplus\Omega^1_{\PP^2}\oplus\OO_{\PP^2}(-3).$$

\subsection{Covering maps}

A covering map of degree $d$ is a flat and finite morphism between schemes, $\varphi\colon X\rightarrow Y$, such that $\varphi_*\OO_X$ is a sheaf of locally free $\OO_Y$-modules with rank $d$.

It is a well known result that for a covering map $\varphi_*\OO_X$ splits as $\OO_Y\oplus\EE$, where $\EE$ is the submodule of trace-zero elements, or Tschirnhausen module, of $\varphi$. Let us recall how this decomposition is obtained, as we will need it for our constructions.

For a covering map of degree $d$, for all $x\in \OO_X$, define $M_x$ as the matrix corresponding to the multiplication by $x\colon \varphi_*\OO_X\rightarrow\varphi_*\OO_X$ for a given basis of $\varphi_*\OO_X$. Then we have a morphism $\tr\colon\varphi_*\OO_X\rightarrow\OO_Y$ defined as $x\mapsto (\operatorname{trace} M_x)$. The map is well defined as the trace of a matrix is unchanged by a change of variables. Furthermore, this morphism splits as for all $y\in \OO_Y$, $\tr(\operatorname{diag}(y,y,\dots,y))=y$. Hence we have the decomposition $\varphi_*\OO_X=\OO_Y\oplus\EE$, where $\EE$ is the kernel of $\tr$.

In the seminal paper \cite{triple}, the case of triple covers is extensively studied. A triple covering map such that $\varphi_*\OO_X=\OO_Y\oplus\EE$ is determined by a morphism $\Phi\colon\Sym^2\EE\rightarrow\varphi_*\OO_X$ that locally can be written as 
\begin{equation}
    \begin{array}{rcr}
        \Phi(y^2) & = & ay+bz+2(a^2-bd)\\
        \Phi(yz) & = & -dy-az-(ad-bc)\\
        \Phi(z^2) & = & cy+dz+2(d^2-ac),
    \end{array}
\end{equation}
where $\{y,z\}$ is a local basis of $\EE$, and $a,b,c,d\in\OO_Y$. Notice that $\Phi$ is determined by $\Phi_2\colon\Sym^2\EE\rightarrow\EE$. The converse is also true. A homomorphism $\Phi_2\in\Hom(\Sym^2\EE,\EE)$ that is locally of the form
\begin{equation}\label{Phi2}
    \begin{array}{ccc}  \Phi_2(y^2) = ay+bz, & \Phi_2(yz) = -dy-az, & \Phi_2(z^2) = cy+dz
\end{array}\end{equation}
induces a commutative and associative multiplication on $\OO_X=\OO_Y\oplus\EE$ whose Tschirnhausen module is $\EE$. Miranda defines a \textit{triple cover homomorphism} as a morphism $\Phi_2\in\Hom(\Sym^2\EE,\EE)$ which, for a local basis of $\EE$, can be written in the form (\ref{Phi2}).

Following Miranda's approach, in \cite{CoverHomomorphisms} a generalisation of this concept is introduced. 
Let $\FF$ and $\Q\subset\RR(\FF)=\bigoplus_n\Sym^n\FF$ be locally free $\OO_Y$-modules with constant rank, $\{z_1,\dots,z_r\}$ and $\{q_i(z_1,\dots,z_r)\}_{i\in I}$ local basis for $\FF$ and $\Q$, respectively. A homomorphism 
$$\Phi\in \Hom\left(\Q,\RR(\FF)/\langle\Q\rangle\right)$$
defines and is defined by a covering map if, locally, the ideal
$$(q_i-\Phi(q_i))\subset k[z_1,\dots,z_r]$$
is a flat deformation of the ideal $(q_i)$. 

If $(q_i)$ is generated by polynomials of degree $2$, i.e. $\Q\subset\Sym^2\FF$, then the (local) equations are 
$$\Phi(q_i)=\sum c_{ij}z_j+d_i,\,\,\, c_{ij},d_i\in\OO_Y.$$

Furthermore, if we assume that $(q_i)$ has a linear first syzygy matrix, then one can prove that the $d_i$ can be written as quadratic forms in the $c_{ij}$, and the relations that the $c_{ij}$ have to satisfy for it to define a covering space are of degree at most $2$. Denote this relations by $I_q\subset k[c_{ij}]$. A homomorphism $\Phi_2\in\Hom(\Q,\FF)$ is called a covering homomorphism if for all open sets, and a choice of basis, it is of the form 
$$\Phi_2(q_i)=\sum_{j=1}^rc_{ij}z_j,$$
with the $c_{ij}$ satisfying the relations in $I_q$.

In the last Section of \cite{CoverHomomorphisms}, the case of a degree $6$ cover $\varphi\colon X\rightarrow Y$, such that 
$$\varphi_*\OO_X=\OO_Y\oplus M\oplus M\oplus \bigwedge^2 M,$$
where $M$ is a rank $2$, simple $\OO_Y$-module, is studied. This is the case we are interested in.

To be more specific, notice first that, as $\OO_{\PP^2}(1)$ is a line bundle over $\PP^2$ for which 
$$\omega_{\PP^2} = \OO_{\PP^2}(1)^{\otimes m}, \omega_{S}\cong \varphi^*\OO_{\PP^2}(1)^{\otimes n},$$ for integers $m,n$, the covering map is Gorenstein, and there is an embedding $$S\hookrightarrow \SSpec \Sym^\bullet \left(\Omega^1\oplus\Omega^1\right)^\vee.$$ 
Moreover, as $\Omega^1_{\PP^2}$ is a simple $\OO_{\PP^2}$-module, $\varphi$ determines and is determined by a \textit{cover homomorphism} in $\Hom\left(\left(\Sym^2\Omega^1_{\PP_2}\right)^{\oplus 3},\left(\Omega^1_{\PP_2}\right)^{\oplus 2}\right)$ and a general fibre of $\varphi$ is given by the polynomials defining the spinor embedding of the affine orthogonal Grassmann variety $\oGr(5,10)$ in $\PP^{15}$. These results appear in Theorem \ref{abelianstructure} with a sketch of the proof. 

We start by describing the vector bundle 
$\SSpec \Sym^\bullet \left(\Omega^1\oplus\Omega^1\right)^\vee$ as a graded ring in Section \ref{Preliminaries}. This will be our starting point. Having chosen coordinates for the ambient space, in Section \ref{linearalgebra} we study the structure of the morphisms $\Hom\left(\Sym^2\Omega^1_{\PP^2},\Omega^1_{\PP^2}\right)$. These constructions are enough to compute the ring $\RR(S,\LL)$ but, as we will see, the number of coefficients and relations between them is still quite large. In order to work them out we will use that, assuming that $\LL$ is symmetric, the coordinate ring of a $(1,3)$-abelian surface is invariant by the action of the extended Heisenberg group $H(1,3)^e$.

\subsection{Canonical level structure} Following \cite{BLangebook,Mumford1}, we present some of the basic theory of abelian varieties needed for our work. As throughout the rest of the paper, let $S$ be an abelian surface and $\LL$ be an ample line bundle defining a polarisation of type $(1,3)$.

Consider the map $\phi_\LL\colon S\rightarrow \hat S$ given by $x\mapsto t_x^*\LL\otimes \LL^{-1}$, where $t_x\colon S\rightarrow S$ is the translation by $x\in S$. Then its kernel $K(\LL)$ is isomorphic to $\ZZ/3\ZZ\times\mu_3$.
The \textit{theta group} $G(\LL)$ is the set of pairs $(x,\phi)$, where $x\in K(\LL)$ and $\phi$ is an isomorphism $\phi\colon \LL\rightarrow t_x^*\LL$. It is a group with multiplication given by
$(y,\psi)\circ (x,\phi)=(x+y,t^*_x\psi\circ\phi)$
and there is an exact sequence 
$$
1\rightarrow k^*\rightarrow G(\LL)\rightarrow K(\LL)\rightarrow 0.
$$
The elements in the kernel are the isomorphisms of $\LL$ with itself.

For each $x\in K(\LL)$ there is an isomorphism $\phi\colon\LL\xrightarrow{\sim} t_x^*\LL$, so we have a representation $\rho\colon K(\LL)\rightarrow\PGL\left(H^0(S,\LL)\right)$. This representation lifts to a linear representation of $G(\LL)$,
\begin{equation*}
\xymatrix{
1\ar[r]  & k^* \ar[r] \ar@{=}[d] & G(\LL) \ar[r]  \ar[d]^{\widetilde\rho} & K(\LL) \ar[r] \ar[d]^{\rho} & 0 \\
1\ar[r] & k^* \ar[r] & \GL(H^0(\LL)) \ar[r] & \operatorname{PGL}(H^0(\LL)) \ar[r] & 0. 
}
\end{equation*}

As an abstract group, $G(\LL)$ is isomorphic to the Heisenberg group $H(1,3)$, and the representation $\rho$ is isomorphic to the Schr\"odinger representation of $H(1,3)$ on $V(\ZZ/3\ZZ)$, the vector space of $k$-valued functions on $\ZZ/3\ZZ$. 

An isomorphism $\theta\colon G(\LL)\xrightarrow{\sim}H(1,3)$ that restricts to the identity on $k^*$ is called a \textit{theta structure}. If it induces a symplectic isomorphism between $K(\LL)$ and $\ZZ/3\ZZ\times\mu_3$, then it is called a \textit{level structure of canonical type on $(S,\LL)$}.

Since the map $\varphi\colon S\rightarrow \PP^2$ is independent of the group law on $S$, by translation, we may assume that $\LL$ is symmetric, i.e. $\LL\cong (-1)^*_S\LL$. Then, $(-1)_S$ induces an automorphism on $K(\LL)$. Denote $K(\LL)^e = K(\LL)\rtimes\langle (-1)_S\rangle$, and define the extended theta group $G(\LL)^e$, as the extension of $K(\LL)^e$ by $k^*$. The extended Heisenberg group is $H^e(1,3)=H(1,3)\rtimes \langle \iota\rangle$, where $\iota$ is the corresponding involution.

As we are considering $\LL$ to be symmetric, we can extend the theta structure to $\theta\colon G^e(\LL)\xrightarrow{\sim}H^e(1,3)$. It induces an isomorphism between $K^e(\LL)$ and $\ZZ/3\ZZ\times\mu_3\times \ZZ/2\ZZ$.

Let $\sigma', \tau'$ and $\iota'$ be elements in $G^e(\LL)$ such that their image are generators of $K^e(\LL)$. By abuse of notation we will denote $\rho(\sigma'),\rho(\tau')$ and $\rho(\iota')$ simply by $\sigma',\tau'$ and $\iota'$.
We can then choose coordinates $\{x_0,x_1,x_2\}$ in $\PP(H^0(S,\LL)^\vee)$, the space of hyperplanes in $H^0(S,\LL)$, such that 
$$\sigma'(x_i)=x_{i+1},\,\,\, \tau'(x_i)=\xi^i x_i,\,\,\, \iota'(x_i)=x_{-i},$$
where $\xi$ is a primitive cubic root of unity, and the indices are taken modulo $3$.
 
Notice that if $(x,\phi)$ is an element of $G^e(\LL)$, then $(x,\phi^n)$ is an element of $G^e(\LL^n)$, so it is natural to extend the action of $H^e(1,3)$ to the ring $\RR(S,\LL)$. To do so, in Section \ref{groupactions} we study the homomorphisms in $\Hom\left(\Sym^2\Omega^1_{\PP^2},\Omega^1_{\PP^2}\right)$ that are equivariant by the action of $H^e(1,3)$.

We then find that the embeddings of $S$ into $\AA\left((\Omega^1_{\PP^2}\oplus\Omega^1_{\PP^2})^\vee\right)$
are related by the group of \mbox{automorphisms} of three distinct points in a plane, isomorphic to $S_3$. This proves our second result.
There is an open and dense subset of the moduli space of abelian surfaces polarised by an ample symmetric line bundle of type $(1,3)$, with canonical level structure, isomorphic to $k^3//S_3\cong k^3$. In other words, $\mathcal{A}^{lev}_{(1,3)}$ is rational.

To finalise we will say some words about how to adapt this construction to the description of irregular surfaces with invariants $p_g(S)=3, q(S)=1, K_S^2=6$.

%% file: RelativeAffineSpacePn.tex
Let $\mathcal{F}$ be a locally free $\OO_{\PP^n}$-module of rank $r$ and $\mathcal{P}_\bullet$ a (fixed) free presentation of $\FF$
\begin{equation}\label{prf} 
0\leftarrow \mathcal{F}\leftarrow \mathcal{P}_0\xleftarrow{M} \mathcal{P}_1 
\end{equation}
where $\mathcal{P}_i=\bigoplus_{j=1}^{m_i} \OO_{\PP^n}(-a_{ij})$ and $M$ is a $m_0\times m_1$ matrix with homogeneous entries in $k[x_0,\dots,x_n]$. Assume that for all $i,j$, $a_{i,j}\geq 0$ and, for simplicity sake, denote $m_0$ by $m$ and $a_{0j}$ by $a_j$.

Denoting by $(x_0,\dots,x_n)$, $(\xi_1,\dots,\xi_m)$ and $\lambda$ a choice of coordinates for $\AA^{n+1}$, $\AA^{m}$ and $\mathbb{G}_m$, we define $\AA(\mathcal{F}^\vee)$ as the quotient of $\mathcal{Z}_M$ by an action of $\mathbb{G}_m$, where
$$\mathcal{Z}_M:=\{\bar\xi M=0\}\subset (\AA^{n+1}\backslash 0)\times \AA^{m}  \text{ , } \bar\xi=(\xi_0,\dots,\xi_m), $$
and the action of $\lambda\in\mathbb{G}_m$ is given by
\begin{equation}\label{action}
\lambda\colon(x_0,\dots,x_n;\xi_0,\dots,\xi_m)\mapsto (\lambda x_0,\dots,\lambda x_n;\lambda^{a_1}\xi_1,\dots,\lambda^{a_m}\xi_m).
\end{equation}
The action is well defined as all polynomials in $\bar\xi M$ are homogeneous in $\lambda$. Furthermore, the ratio $(x_0:x_1:\dots:x_n)$ is preserved by the action of $\mathbb{G}_m$ so the projection to the first factor defines a morphism $\pi\colon\AA\left(\FF^\vee\right)\rightarrow\PP^n$.
$$
\xymatrix{
(\AA^{n+1}\backslash 0)\times \AA^{m} \ar[d]_{p_1}  \ar[r] & \AA\left(\FF^\vee\right) \ar[d]_\pi \\
(\AA^{n+1}\backslash 0) \ar[r] & \PP^n
}$$

\begin{definition}\label{projectiveF}
 Given a locally free $\OO_{\PP^n}$-module $\FF$ of rank $r$, and $\mathcal{P}_\bullet$ a presentation of $\FF$, 
$$ \mathcal{F}\leftarrow \bigoplus_{i=1}^m\OO_{\PP^n}(-a_i)\xleftarrow{M} \bigoplus_{i=1}^{m_1}\OO_{\PP^n}(-a_{1i})$$
we denote by $R\left(\mathcal{P}_\bullet\right)$ the ring $k[x_0,\dots,x_n;\xi_1,\dots,\xi_m]/(\overline{\xi}M)$ and denote the scheme $\Proj\left(R(\mathcal{P}_\bullet)\right)$, where we consider the irrelevant ideal to be $m=(x_0,\dots,x_n)$, as
$$\pi\colon\AA\left(\FF^\vee\right)\rightarrow\PP^n.$$
\end{definition}

Notice that $\FF$ is the cokernel of the map $M\colon\mathcal{P}_1\rightarrow\mathcal{P}_0$ and hence, any morphism in $\Hom\left(\FF,\OO_{\PP^n}\right)$ is given by a morphism $\overline\xi\in\Hom\left(\mathcal{P}_0,\OO_{\PP^n}\right)$ for which $\overline{\xi}M=0$.

It is important to note that we defined $\AA\left(\FF^\vee\right)$ for a locally free sheaf $\mathcal{F}$, but for the construction we use a free resolution which is not unique. As an example take $\mathcal{F}=\OO_{\PP^1}(-a)$, $a\geq 0$, for which we can define $\PP\left(\FF^\vee\right)$ using both resolutions
$$\begin{array}{rcl}
\mathcal{F} & \leftarrow & \OO_{\PP^1}(-a) \leftarrow 0, \\
\mathcal{F} & \leftarrow & \OO_{\PP^1}(-a-1)^{\oplus 2}\leftarrow \OO_{\PP^1}(-a-2) \leftarrow 0.
\end{array}$$

As for the example above, as a scheme, $\PP\left(\FF^\vee\right)$ is independent of the resolution for $\FF$. This is stated in the next theorem.

\begin{theorem}\label{theoremproj}Let $\mathcal{P}_\bullet$, $\mathcal{Q}_\bullet$ be free resolutions of $\mathcal{F}$, a locally free $\OO_{\PP^n}$-module of rank $r$. Then
 \begin{enumerate}
  \item \label{unprojthm}
  As schemes over $\PP^n$, 
 $$\big(\pi\colon\AA(\mathcal{P}_\bullet)\rightarrow\PP^n\big)\cong \big(\pi'\colon\AA(\mathcal{Q}_\bullet)\rightarrow\PP^n\big).$$
 \item \label{ring-manifold}
There is a set of elements $\gamma_i\in k[x_i]$ such that for any $i\in I$, $I$ a index set, we have 
$$R(\mathcal{P}_\bullet)_{\gamma_i}\cong k[\gamma_i^{-1},x_0,\dots,x_n;\xi_{(1,i)},\xi_{(1,i)}\dots,\xi_{(r,i)}]$$
where each $\xi_{(j,i)}$ is a $\OO_{\PP^n}$-linear combination of the $\xi_j$.
 \end{enumerate}
\end{theorem}

\begin{proof}
$(1)$ Assume $\mathcal{P}_\bullet$ is the minimal free resolution of $\mathcal{F}$. Then we have a morphism of chain complexes,
$$\xymatrix{ 
\mathcal{F} \ar@{=}[d] & \ar[l] \mathcal{Q}_0 \ar[d]_{f_0} & \ar[l]_{M_{\mathcal{Q}}} \mathcal{Q}_1 \ar[d]_{f_1} &  \\
\mathcal{F} & \ar[l] \mathcal{P}_0 & \ar[l]_{M_{\mathcal{P}}} \mathcal{P}_1 &  \\
}$$
that induces a morphism between the bi-graded rings  $R(\mathcal{Q}_\bullet)$ and $R(\mathcal{P}_\bullet)$ as, by minimality of the free resolution $\mathcal{P}_\bullet$, any morphism in $\Hom(\mathcal{Q}_0,\OO_{\PP^n})$ can be factored as $f_0$ and a morphism in $\Hom(\mathcal{P}_0,\OO_{\PP^n})$ that we denote by $f(\xi)$ for all $\xi\in\Hom(\mathcal{Q}_0,\OO_{\PP^n})$. Fixing a basis of $\OO_{\PP^n}$ we have a map
$$\begin{matrix}
f\colon & R(\mathcal{Q}_\bullet) & \longrightarrow & R(\mathcal{P}_\bullet) \\
 & x_i & \mapsto & x_i\\
 & \xi_i & \mapsto & f(\xi_i), 
\end{matrix}$$
where $\xi_i$ is a basis of $\Hom(\mathcal{P}_0,\OO_{\PP^n})$, hence a morphism $\AA(\mathcal{Q}_\bullet)\rightarrow\AA(\mathcal{P}_\bullet)$.

As $\FF$ is locally free, take $\mathcal{U}\subset\PP^n$ an open set such that 
$$\AA(\mathcal{P}_\bullet)|_{\mathcal{U}}\cong \Spec\big(\OO_{\PP^n}(\mathcal{U})[\xi_{(1)},\dots,\xi_{(r)}]\big).$$
Then these $\xi_i|_\mathcal{U}$ are a basis for $\FF^\vee|_\UU$ and, by the commutativity of the first square in the diagram, $f_0|_\mathcal{U}$ is an isomorphism. We get that for any two resolutions of $\FF$, $\mathcal{Q}^1_\bullet$, $\mathcal{Q}^2_\bullet$, we have isomorphisms
$$\AA(\mathcal{Q}^1_\bullet)\xleftarrow{d_1} \AA(\mathcal{P}_\bullet)\xrightarrow{d_2} \AA(\mathcal{Q}^2_\bullet),$$ 
where $\mathcal{P}_\bullet$ is a minimal free resolution so, even without a morphism in between them, they are isomorphic through $\AA(\mathcal{P}_\bullet)$.

For $(2)$ just pick $\gamma_i$ general enough and apply Gauss elimination to $M_{\mathcal{P}}$. 
\end{proof}

With this theorem we can indeed write $\AA(\FF^\vee)$ and, by (\ref{ring-manifold}), we see that the construction is $\SSpec \Sym^\bullet \FF^\vee$.  Most of the time we will be using the construction with respect to the minimal free resolution of $\mathcal{F}$. When not, we will be working locally so there will be no harm done by part (\ref{ring-manifold}) of Theorem \ref{theoremproj}. 

\begin{example}\label{babycaseexample}
Let $\EE=\Omega^1_{\PP^2}(-m)$, $m\geq 0$. Using the Euler sequence
$$\Omega^1(-m)\leftarrow \OO_{\PP^2}(-m-2)^{\oplus 3}\xleftarrow{M} \OO_{\PP^2}(-m-3)\leftarrow 0,$$
where $M=\left(\begin{matrix}x_0 & x_1 & x_2\end{matrix}\right)^t$, we have
$$\AA\left(\Omega^1_{\PP^2}(-m)^\vee\right)=\Proj\left(k[x_0,x_1,x_2;y_0,y_1,y_2]/\left(\sum x_iy_i\right)\right),$$
with $\deg x_i=1$, $\deg y_i=m+2$. Let $\UU_i=\{x_i\neq 0\}$, then over $\UU_2$ we have
$$\AA\left(\Omega^1_{\PP^2}(-m)^\vee\right)\big|_{\pi^{-1}(\UU_2)}=\Proj\left(k[x_0,x_1,x_2^{\pm};y_0,y_1]\right)$$
and the transition map from $\UU_2$ to $\UU_1$ is given by the change of coordinates
$$\left(\begin{matrix}
   y_0   &   y_1  
\end{matrix}\right)\mapsto
\left(\begin{matrix}
 y_0 &  -\tfrac{x_0}{x_2}y_0-\tfrac{x_1}{x_2}y_1
\end{matrix}\right)=
\left(\begin{matrix}
   y_0   &   y_2  
\end{matrix}\right).
$$
\end{example}

\begin{example}\label{unprojexample}
In this example we want to see what the following inclusion 
$$\OO_{\PP^2}(-3)\subset \bigoplus_{n\geq 0}\Sym^n\left(\Omega^1_{\PP^2}\oplus\Omega^1_{\PP^2}\right)$$
implies in terms of the graded ring defining $\AA\left(\left(\Omega^1_{\PP^2}\oplus\Omega^1_{\PP^2}\right)^\vee\right)$. Recall that the inclusion above is expected
\begin{equation*}
\begin{array}{rcl}
\Sym^2\left(\Omega^1_{\PP^2}\right)^{\oplus 2} &=& \Sym^2\left(\Omega^1_{\PP^2}\right)\oplus\Omega^1_{\PP^2}\otimes \Omega^1_{\PP^2}\oplus \Sym^2\left(\Omega^1_{\PP^2}\right)^{\oplus 2} \\
& =& \left(\Sym^2\left(\Omega^1_{\PP^2}\right)\right)^{\oplus 3}\oplus\left(\Omega^1_{\PP^2}\bigwedge \Omega^1_{\PP^2}\right) \\
& =& \left(\Sym^2\left(\Omega^1_{\PP^2}\right)\right)^{\oplus 3}\oplus\OO_{\PP^2}(-3).
\end{array}
\end{equation*}
As for the Example \ref{babycaseexample}, we have the following equality
$$\AA\left(\left(\Omega^1_{\PP^2}\oplus\Omega^1_{\PP^2}\right)^\vee\right)=\Proj\left(k[x_i;y_i,z_i]/\left(\sum x_iy_i,\sum x_iz_i\right)\right)_{0\leq i\leq 2},$$
with $\deg(x_i)=1$, $\deg(y_i)=\deg(z_i)=2$. The interesting fact is that the summand $\OO_{\PP^2}(-3)$ should correspond to a variable of degree three and at the same time to the multiplication of variables of degree two.

The first explanation for this came from using unprojection methods (see \cite{unproj1,unproj2}), i.e. the unprojection of the following pair of ideals
$$\left(\sum x_iy_i,\sum x_iz_i\right)\subset (x_i).$$ 
This example is described in \cite[$\S 4$]{unproj1} but, for completeness sake, we recall the result. 
From the equation
$$\left(\begin{matrix}
y_0 & y_1 & y_2 \\
z_0 & z_1 & z_2
\end{matrix}\right)
\left(\begin{matrix}
x_0 \\ x_1 \\ x_2
\end{matrix}\right)=0,
$$
we know that the first matrix has rank $1$ and, by Cramer's rule, also the following vector is in its kernel  
$$\bigwedge^2
\left(\begin{matrix}
y_0 & y_1 & y_2 \\
z_0 & z_1 & z_2
\end{matrix}\right)=
\left(\begin{matrix}
y_1z_2-y_2z_1 \\ y_2z_0-y_0z_2 \\ y_0z_1-y_1z_0
\end{matrix}\right).$$
We conclude that there exists an element $t$ such that
$$\bigwedge^2
\left(\begin{matrix}
y_0 & y_1 & y_2 \\
z_0 & z_1 & z_2
\end{matrix}\right)
=t
\left(\begin{matrix}
x_0 \\ x_1 \\ x_2
\end{matrix}\right).
$$
The degree of $t$ is $3$, and this is the quadratic term corresponding to the summand $\OO_{\PP^2}(-3)$.

A different reasoning comes from applying Theorem \ref{theoremproj} to the following free resolutions
\begin{equation*}
\xymatrix{ 
\Omega^1_{\PP^2}\bigwedge\Omega^1_{\PP^2} \ar@{=}[d]  & \ar[l] \OO_{\PP^2}(-4)^{\oplus 3} \ar[d]_{f_0} & \ar[l] \OO_{\PP^2}(-5)^{\oplus 3}  &   \\
\OO_{\PP^n}(-3) \ar@{=}[r] & \OO_{\PP^n}(-3)  & &    \\
}
\end{equation*}
The map $f_0$ is exactly the map $\left(\begin{matrix}x_0 & x_1 & x_2\end{matrix}\right)$. 
\end{example}

We included the relation between the unprojection method and the description of the ring $R\left(\Omega^1_{\PP^2}\oplus\Omega^1_{\PP^2} \right)$ as we believe that describing similar rings for higher dimensions, e.g. $R\left(\left(\Omega^1_{\PP^3}\right)^{\oplus 3} \right)$, besides the explicit construction of an ambient space for (some) irregular $3$--folds, might be the source of new unprojection methods, see \cite{PapadakisII,PapadakisIII,PapadakisGen}.

%% file: CoverHomomorphism.tex
The goal of this section is to study trace-zero homomorphisms in $\Hom(\Sym^2\Omega^1_{\PP^2},\Omega^1_{\PP^2})$ that are equivariant as $H(1,3)^e$-sheaves of $\OO_{\PP^2}$-modules. We present most of the results for the twisted sheaves $\Omega^1_{\PP^2}(-m)$ so they can be used to describe projective models of irregular surfaces, Section \ref{irregular}.

\subsection{$\Hom\left(\Sym^2\Omega^1_{\PP^2},\OO_{\PP^2}\right)$}

Following Example \ref{babycaseexample},
let $(y_0,y_1,y_2)$ be a basis for $\Omega^1_{\PP^2}(-m)^\vee$ and write a basis for $\Sym^2\left(\Omega^1_{\PP^2}(-m)\right)^\vee$ as the symmetric matrix 
$$Y^2:=\left(\begin{matrix}
  y_0^2 & y_0y_1 & y_0y_2 \\
  y_0y_1 & y_1^2 & y_1y_2 \\
  y_0y_2 & y_1y_2 & y_2^2
 \end{matrix}\right).$$
The relation given by $M$ translates as $Y^2\overline{x}=0$, where $\overline{x}=\left(\begin{matrix}x_0 & x_1 & x_2 \end{matrix}\right)^t$.
So we need to find a $3\times 3$ symmetric matrix, $L$, such that $L\overline{x}=0$.

\begin{lemma}[Kova\v{c}ec Lemma]\label{kovaceclemma}
If $L$ is a $3\times 3$ symmetric matrix in $k[x_0,x_1,x_2]$ such that $L \overline{x}=0$ then $L=M_3 N M_3$ where $N$ is a symmetric matrix and $M_3$ is the syzygy matrix of the ideal $(x_0,x_1,x_2)$,
$$M_3=\left(\begin{matrix}
    0 & x_2 & -x_1 \\
    -x_2 & 0 & x_0 \\
    x_1 & -x_0 & 0
   \end{matrix}\right).
$$
\end{lemma}

\begin{proof}It follows from the Koszul complex that $L \overline{x}=0$ implies the existence of a $3\times3$ matrix $P$ such that $L=P M_3$, therefore
\begin{equation*}
\begin{aligned}
 0 & =  L \overline{x}  =  L^t \overline{x}  =  -M_3 P^t \overline{x} & \implies & P^t \overline{x}=\lambda \overline{x}, & &  \\
  &&\implies & (P^t-\lambda I)=Q M_3 \\
  &&\implies & L = -M_3 Q M_3 - \lambda M_3 \\
  &&\iff & L^t = -M_3 Q^t M_3 + \lambda M_3 , && \lambda\in k[x_0,x_1,x_2].
 \end{aligned}
\end{equation*}
Taking $N=-\frac{(Q+Q^t)}{2}$ we get the result.
\end{proof}

\begin{corollary}\label{zerohom}
A morphism $\Phi\in \Hom(\Sym^2(\Omega^1_{\PP^2}(-m_1)),\OO_{\PP^2}(-m_2))$ can be written as a $3\times 3$ matrix,
$$\Phi\left(Y^2\right)=M_3\left(
\begin{matrix}
 d_1 & d_2 & d_3 \\
     & d_4 & d_5 \\
     &     & d_6
\end{matrix}\right)M_3$$
where $d_i\in k[x_0,x_1,x_2]$, $\deg(d_i)=2m_1-m_2+2$. In particular, $\Hom(\Sym^2\Omega^1_{\PP^2},\OO_{\PP^2}(-3))=0$. 
\end{corollary}

\subsection{$\Hom\left(\Sym^2\Omega^1_{\PP^2},\Omega^1_{\PP^2}\right)$}\label{triplecoverhomo}

Using Lemma \ref{kovaceclemma} we now can describe $\Hom(\Sym^2(\Omega^1_{\PP^2}(-m)),\Omega^1_{\PP^2}(-m))$ and the change of basis so that it is defined by a trace free basis. 

Let $\Phi \in\Hom\left(\Sym^2(\Omega^1_{\PP^2}(-m)),\Omega^1_{\PP^2}(-m)\right)$ and $\widetilde{\Phi}\in\Hom\left(\OO_{\PP^2}(-2m-4)^{\oplus 6},\Omega^1_{\PP^2}(-m)\right)$ such that $\widetilde{\Phi}\circ M=0$.  
\begin{equation*}
\xymatrix{
\Sym^2(\Omega^1_{\PP^2}(-m)) \ar[d]_\Phi & \ar[dl]_{\widetilde{\Phi}} \ar[d]_{\Phi'} \ar[l]_p \OO_{\PP^2}(-2m-4)^{\oplus 6} & \ar[l]_M \OO_{\PP^2}(-2m-5)^{\oplus 3} & \ar[l] 0 \\
\Omega^1_{\PP^2}(-m) & \ar[l]_q \OO_{\PP^2}(-m-2)^{\oplus 3} & \ar[l]_N \OO_{\PP^2}(-m-3) & \ar[l] 0
}
\end{equation*}
Applying the functor $\Hom\left(\OO_{\PP^2}(-2m-4)^{\oplus 6},-\right)$ to the lower exact sequence we get 
$$0=H^1\left(\OO_{\PP^2}(m+1)^{\oplus 6}\right)\leftarrow  H^0\left(\Omega^1_{\PP^2}(m+4)^{\oplus 6}\right) \leftarrow H^0\left(\OO_{\PP^2}(m+2)^{\oplus 18}\right) \leftarrow \dots$$
so for any $\widetilde\Phi$ there exists $\Phi' \in\Hom( \OO_{\PP^2}(-2m-4)^{\oplus 6},\OO_{\PP^2}(-m-2)^{\oplus 3})$ such that $\widetilde\Phi=q\circ \Phi'$. Given a basis $(x_i,y_i)$ for $\PP^2$, $\Omega^1_{\PP^2}(-m)^\vee$, we write $\Phi'$ as
$$\Phi'\left(Y^2\right)=C_0y_0+C_1y_1+C_2y_2$$ 
where $C_0, C_1, C_2$ are symmetric matrices in $k[x_i]$. The ideal case would be to have $C_i=M_3N_iM_3$ but this is not necessarily true, just consider $C_i=x_iC$ for some symmetric matrix $C$. On the other hand for this case we have $Y^2\mapsto C\left(\sum x_iy_i\right)=0$ and this comes from the fact that $\widetilde\Phi\circ M=0\Leftrightarrow q\circ \Phi'\circ M=0\Rightarrow \Phi'\circ M\in \operatorname{Im} N$.

Although we have this ambiguity in global terms, if $\Phi_1'\left(Y^2\right)=\sum C_iy_i$ and $\Phi_2'\left(Y^2\right)=\sum C'_iy_i$ differ by an element $\sum (x_iC)y_i\in \operatorname{Im} N$, then locally they are the same. E.g. take the open set $\UU=\{x_2\neq 0\}$ then
$$\Phi'_1\left(Y^2\right)|_\UU=\sum_{i=0}^{1} \left(C_i-\frac{x_i}{x_2}C_2\right)y_i=\sum_{i=0}^{1} \left(C'_i-\frac{x_i}{x_2}C'_2\right)y_i=\Phi_2'\left(Y^2\right)|_\UU$$
and the morphisms to any other open set are just given by an $\OO_{\PP^2}$-automorphism.

\begin{proposition}\label{proposition22}
 Given a morphism $\Phi\in\Hom\left(\Sym^2(\Omega^{1}_{\PP^2}(-m)),\Omega^{1}_{\PP^2}(-m)\right)$ and $(x_i)_{0\leq i\leq 2}$ a basis for $\OO_{\PP^2}$. Then there is a basis $(y_i)_{0\leq i\leq 2}$ for $\Omega^{1}_{\PP^n}$ such that $\Phi$ can be written as
 $$\Phi\left(Y^2\right)= C_0y_0+C_1y_1+C_2y_2,$$
 where $C_i\overline{x}=0$ for all $i$.
\end{proposition}

\begin{proof}Assume $\Phi$ can be written as $C'_0y_0+C'_1y'_1+C'_2y_2$ for a basis $(x_i,y_i)_{0\leq i\leq 2}$ of $\PP^2$ and $\Omega^{1}_{\PP^2}$. If $C'_2\overline{x}=0$ then 
$$\Phi|_{\{x_2\neq 0\}}\circ M=0\Leftrightarrow (x_2C'_i-x_iC'_2)\overline{x}=0\Rightarrow C'_i\overline{x}=0,$$
so we just need to make a change of variables in the $y_i$ such that $C'_2\overline{x}=0$.

Decomposing $C'_2$ as $C_2+x_2C$, where the entries of $C_2$ are in $k[x_0,x_1]$ for $i\in\{0,1\}$, define $C_i=C'_i-x_iC$.
Let $\widehat{C}_{2}$ be the submatrix of $C_2$ with its first $2$ columns and $2$ rows. Then $\widehat{C}_{2}$ satisfies $\widehat{C}_{2}\left(\begin{smallmatrix}x_0 & x_1\end{smallmatrix}\right)^t=0$ as none of its entries has a multiple of $x_2$ and $(C_i-x_i/x_2C_2)\overline{x}=0$. 

We only have to deal with the last column and row of $C_2$. Choosing a decomposition of $C_2{[3,3]}=2\left(a_0x_0+a_1x_1\right)$, $a_i$ terms of degree one in $k[x_0,x_1]$, and changing variables $y_2\mapsto y_2-(a_0x_0+a_1x_1)$ and $y_i\mapsto y_i+a_ix_2$ for $i \neq 2$, only the last row and column of $C_2$ change and $C_2[3,3]=0$. $C_2$ has the form
$$\left(\begin{matrix}
   \widehat{C}_{2} & v^t \\
   v & 0
  \end{matrix}\right)$$
where $\widehat{C}_{2}\left(\begin{smallmatrix}x_0 & x_1\end{smallmatrix}\right)^t=0$. Using the same trick as in the beginning, we can get rid of all the multiples of $x_2$ in $v$. With all the entries of $v=\left(\begin{smallmatrix}v_0 & v_1\end{smallmatrix}\right)$ in $k[x_0,x_1]$ and $v$ satisfying $v_0x_0+v_1x_1=0$, a change of variables $y_i\mapsto y_i-v_i$ completes the proof.
\end{proof}

In the next proposition we show how to change variables so that the basis of $\Omega^1_{\PP^2}(-m)^\vee$ we are using is trace-free. Notice that being trace-free is a local property but the form $\Phi$ is written above is a global decomposition that might not be trace-free. 

\begin{corollary}\label{tracefree}
Given $\Phi\in\operatorname{Hom}\left(\Sym^2(\Omega^{1}_{\PP^2}(-m)),\Omega^{1}_{\PP^2}(-m)\right)$ and any basis $(x_i)_{0\leq i\leq 2}$ for $\OO_{\PP^2}$, there is a basis $(y_i)$ for $\Omega^{1}_{\PP^2}(-m)^\vee$ such that we can write $\Phi\left(Y^2\right)$ as $\sum C_iy_i$ and for all open sets $\{x_i\neq 0\}$ the variables $y_j$, $j\neq i$, are trace zero elements.
\end{corollary}

\begin{proof}
Using Proposition \ref{proposition22}, write $\Phi\left(Y^2\right)$ as $\sum C_iy_i$ where $C_i\overline{x}=0$ for all $i$. If $x_2\neq 0$, we have that the trace of $y_i$, $i\neq 2$, is given by 
$$\operatorname{tr}(y_i)=\sum_{j\neq 2}\left(C_j-\frac{x_j}{x_2}C_2\right){[i,j]}=\sum_{j=0}^2\left(C_j-\frac{x_j}{x_2}C_2\right){[i,j]}=\sum_{j=0}^2C_{j}[i,j]$$ 
In particular, $\operatorname{tr}(y_i)$ does not depend on the open $\{x_i\neq 0\}$ and we can globally define $v_i=\operatorname{tr}(y_i)$. From the following equalities
$$\sum_i x_iv_i=\sum_ix_i\sum_{j}C_{j}[i,i+j]=\sum_j\sum_ix_iC_j[i,i+j]=0$$
one see that we can use the entries $v_i$ to write a vector $v$ as $v=v'M_n$ which allows to use the change of variables $y_i\mapsto y_i-\frac{v_i}{3}$ and we are done.
\end{proof}

One can then write $\Phi\in\operatorname{Hom}\left(\Sym^2(\Omega^1_{\PP^2}(-m)),\Omega^1_{\PP^2}(-m)\right)$ as 
$$\Phi\left(Y^2\right)=C_0y_0+C_1y_1+C_2y_2,\hspace{3mm}C_i=M_3\left(\begin{matrix}c_{i1} & c_{i2} & c_{i3} \\ & c_{i4} & c_{i5} \\ sym & & c_{i6} \end{matrix}\right)M_3$$
where each $c_{ij}$ is a polynomial of degree $m$. Using the change of variables described in Corollary \ref{tracefree}, $\Phi|_{\{x_2\neq 0\}}$ can be written as 
\begin{equation}\label{equations}
\begin{array}{rcr}
y_0^2 & = & c_1y_0+c_0y_1\\ 
y_0y_1 & = & -c_2y_0-c_1y_1 \\
y_1^2 & = & c_3y_0 + c_2y_1
\end{array}
\end{equation}
where each $c_i$ is of the form 
\begin{equation}
\begin{aligned}
 c_0  = & 
 x_{2}^{2} (c_{14}) 
 -  x_{1} x_{2}( 2c_{15} + c_{24}) 
 + x_{1}^{2}  (c_{16} + 2 c_{25}) 
 -  \frac{x_{1}^{3}}{x_{2}} (c_{26})\\
 c_1  = &
 \tfrac{1}{3} x_{2}^{2}( c_{04} + 2 c_{12}) 
 - \tfrac{2}{3} x_{1} x_{2} (c_{05} + c_{13} + c_{22})
 + \tfrac{1}{3} x_{1}^{2} (c_{06} + 2 c_{23}) \\
 &  - \tfrac{1}{3} x_{0} x_{2}( 2c_{15} + c_{24}) 
 + \tfrac{2}{3} x_{0} x_{1}( c_{16} + 2 c_{25})
 - \frac{x_{0} x_{1}^{2}}{x_{2}}(c_{26}) \\
 c_2  = & 
 \tfrac{1}{3} x_{2}^{2}( 2c_{02}  + c_{11})
 - \tfrac{1}{3} x_{1} x_{2} (2c_{03} + c_{21})
 - \tfrac{2}{3} x_{0} x_{2} (c_{05}  + c_{13} + c_{22}) \\
 & 
 + \tfrac{2}{3} x_{0} x_{1} (c_{06} + 2 c_{23})
 + \tfrac{1}{3} x_{0}^{2} (c_{16} + 2c_{25}) 
 -  \frac{x_{0}^{2} x_{1}}{x_{2}} (c_{26}) \\
 c_3  = & 
  x_{2}^{2} (c_{01}) 
 -  x_{0} x_{2} (2c_{03} + c_{21}) 
 + x_{0}^{2} (c_{06} + 2 c_{23})
 -\frac{  x_{0}^{3}}{x_{2}} (c_{26}).\\
\end{aligned}
\end{equation}

From the open set $\{x_2\neq 0\}$ to any other $\{x_1\neq 0\}$ it is just a $\OO_{\PP^2}$--linear transformation, so the $10$ parameters 
\begin{equation*}
\begin{matrix}
 c_{01}, 2c_{02}  + c_{11}, 2c_{03} + c_{21}, c_{04} + 2 c_{12}, c_{05} + c_{13} + c_{22}, \\
\hspace{20mm} c_{06} + 2 c_{23}, c_{14}, 2c_{15} + c_{24}, c_{16} + 2 c_{25}, c_{26} 
\end{matrix}
\end{equation*}
describe all the covering homomorphisms in $\Hom(\Sym^2(\Omega^1_{\PP^2}(-m)),\Omega^1_{\PP^2}(-m))$. 

To simplify notation we will denote the parameters above as
\begin{table}[h!]
\centering
\begin{tabular}{ll|l|l|ll}
\hline
\multicolumn{1}{|l|}{$c_{01}$} & $\beta_0$ & $2c_{02}  + c_{11}$ & $\beta_{01}$ & \multicolumn{1}{l|}{$c_{05} + c_{13} + c_{22}$} & \multicolumn{1}{l|}{$\beta_{012}$} \\ \hline
\multicolumn{1}{|l|}{$c_{14}$} & $\beta_1$ & $2c_{03} + c_{21}$  & $\beta_{02}$ &                       &                       \\ \cline{1-4}
\multicolumn{1}{|l|}{$c_{26}$} & $\beta_2$ & $c_{04} + 2 c_{12}$ & $\beta_{10}$ &                       &                       \\ \cline{1-4}
                               &           & $2c_{15} + c_{24}$  & $\beta_{12}$ &                       &                       \\ \cline{3-4}
                       &                   & $c_{06} + 2 c_{23}$ & $\beta_{20}$ &                       &                       \\ \cline{3-4}
                       &                   & $c_{16} + 2 c_{25}$ & $\beta_{21}$ &                       &                       \\ \cline{3-4}
\end{tabular}
\end{table}

\noindent so the equations over $\{x_2\neq 0\}$ look like, 
\begin{equation}\label{equationsc}
\begin{aligned}
 c_0  = & 
 \beta_{1}  x_{2}^{2} 
 -  \beta_{12} x_{1} x_{2} 
 +   \beta_{21}  x_{1}^{2}
 -   \beta_{2} \frac{x_{1}^{3}}{x_{2}} \\
 c_1  = &
 \tfrac{1}{3} \beta_{10} x_{2}^{2} 
 - \tfrac{2}{3}  \beta_{012} x_{1} x_{2}
 + \tfrac{1}{3}  \beta_{20} x_{1}^{2} 
   - \tfrac{1}{3} \beta_{12} x_{0} x_{2} 
 + \tfrac{2}{3} \beta_{21} x_{0} x_{1}
 - \beta_{2} \frac{x_{0} x_{1}^{2}}{x_{2}}\\
 c_2  = & 
 \tfrac{1}{3} \beta_{01} x_{2}^{2}
 - \tfrac{1}{3}  \beta_{02} x_{1} x_{2}
 - \tfrac{2}{3}  \beta_{012} x_{0} x_{2}  
 + \tfrac{2}{3}  \beta_{20} x_{0} x_{1}
 + \tfrac{1}{3} \beta_{21} x_{0}^{2}  
 -   \beta_{2} \frac{x_{0}^{2} x_{1}}{x_{2}} \\
 c_3  = & 
  \beta_{0}  x_{2}^{2} 
 -  \beta_{02} x_{0} x_{2}  
 +  \beta_{20} x_{0}^{2}
 - \beta_{2} \frac{x_{0}^{3}}{x_{2}} \\
 \end{aligned}
\end{equation}

\subsection{Group actions}\label{groupactions}

Recall that if $S$ is an abelian surface, $\LL$ is a symmetric line bundle that induces a polarisation of type $(1,3)$, and $\theta$ is a canonical level structure, then the morphism $\varphi_\mathcal{|\LL|}$ is equivariant for the Heisenberg group $H(1,3)^e$. Choosing coordinates $\{x_0,x_1,x_2\}$ for $\PP^2$, an action of $H^e(1,3)$ is generated by 
\begin{equation}\label{groupoperations}
\begin{array}{ccc}
    \sigma'(x_i)=x_{i+1} & \iota'(x_i)=x_{-i} & \tau'(x_i)=\xi^i x_i,   
\end{array}
\end{equation}
where $\xi$ is a primitive cubic root of unity. Notice that $\sigma',\iota'$ generate the symmetric group $S_3$ and $\tau'$ generates the multiplicative group $\mu_3$. In this Section we will extend the action of $H(1,3)^e$ to the graded ring
$$R=k[x_0,x_1,x_2,;y_0,y_1,y_2]/\left(\sum x_iy_i\right)$$ 
with $\deg(x_i)=1,\deg(y_i)=2$, so that the triple covering map $\varphi\colon S'\rightarrow \PP^2$ such that $\varphi_*\OO_{S'}=\OO_{\PP^2}\oplus\Omega^1_{\PP^2}$ is equivariant for $H^e(1,3)$. This will be a fundamental step to finish the description of the ring $\RR(S,\LL,\theta)$.

Consider the following extension
\begin{equation}\begin{array}{rcl}
\sigma(x_i,y_j) & = & (x_{i+1},\xi^{a} y_{j+1}), \\
\iota(x_i,y_j) & = & (x_{-i},(-1)^{b} y_{-j}), \\
\tau(x_i,y_j) & = & (\xi^i x_i,\xi^{c} \xi^{-i}y_i),  
\end{array}\end{equation}
where $a, b, c$ are integers. These are all the extensions which $1)$ have the correct order, $2)$ fix the vanishing of the polynomial $\sum_i x_iy_i$, $3)$ leave the basis $\{y_i\}$ trace-free.

The goal is to prove that we can consider $a=b=c=0$.
To do so, we will examine the action of each of these elements on the coefficients $\beta_i$ of a covering homomorphism  $\Phi_2\in\Hom\left(\Sym^2\Omega^1_{\PP^2},\Omega^1_{\PP^2}\right)$.

Let us start with $\sigma$. Over the open set $\UU_1=\{x_1\neq 0\}$, $\Phi_2|_{\UU_1}$ is of the form,
$$\begin{array}{rcl}
y_0^2 & = & c_1'y_0+c_0'y_2 \\
y_0y_2 & = & -c_2'y_0-c_1'y_2 \\
y_2^2 & = & c_3'y_0+c_2'y_2.
\end{array}$$
Applying $\sigma$ we get 
$$\begin{array}{rcl}
\xi^{2a}y_1^2  & = & \xi^a\sigma(c_1')y_1+\xi^a\sigma(c_0')y_0 \\
\xi^{2a}y_1y_0 & = & -\xi^a\sigma(c_2')y_1-\xi^a\sigma(c_1')y_0 \\
\xi^{2a}y_0^2  & = & \xi^a\sigma(c_3')y_1+\xi^a\sigma(c_2')y_0. 
\end{array}$$
We get the following equalities
$$c_0=\xi^{-a}\sigma(c_3'),\,\, c_1=\xi^{-a}\sigma(c_2'),\,\, c_2=\xi^{-a}\sigma(c_1'),\,\, c_3=\xi^{-a}\sigma(c_0'),$$
where the $c_i$ appear in $\Phi_2|_{\UU_2}$ (see Equation (\ref{equationsc})), and the $c'_i$ in $\Phi_2|_{\UU_1}$ are of the form
\begin{equation}
\begin{aligned}
 c'_0  = & 
 \beta_{2}  x_{1}^{2} 
 +  \beta_{12}  x_{2}^2 
 -   \beta_{21}  x_{1}x_{2}
 -   \beta_{1} \frac{x_{2}^{3}}{x_{1}} \\
 c'_1  = &
 \tfrac{1}{3} \beta_{10} x_{2}^{2} 
 - \tfrac{2}{3}  \beta_{012} x_{1} x_{2}
 + \tfrac{1}{3}  \beta_{20} x_{1}^{2} 
   + \tfrac{2}{3} \beta_{12} x_{0} x_{2} 
 - \tfrac{1}{3} \beta_{21} x_{0} x_{1}
 - \beta_{1} \frac{x_{0} x_{2}^{2}}{x_{1}}\\
 c'_2  = & 
- \tfrac{1}{3} \beta_{01} x_{1}x_{2}
 + \tfrac{1}{3}  \beta_{02} x_{1}^{2}
 + \tfrac{2}{3}  \beta_{10} x_{0} x_{2}  
 - \tfrac{2}{3}  \beta_{012} x_{0} x_{1}
 + \tfrac{1}{3} \beta_{12} x_{0}^{2}  
 -   \beta_{1} \frac{x_{0}^{2} x_{2}}{x_{1}} \\
 c'_3  = & 
  \beta_{0}  x_{1}^{2} 
 -  \beta_{01} x_{0} x_{1}  
 +  \beta_{10} x_{0}^{2}
 - \beta_{1} \frac{x_{0}^{3}}{x_{1}} 
\end{aligned}
\end{equation}
As there is no relation between the $x_i$ one concludes that
\begin{equation}\label{conditionssigma}
\begin{array}{c}
   \beta_0=\xi^{a}\beta_1=\xi^{2a}\beta_2 \\
   \beta_{01}=\xi^{a}\beta_{12}=\xi^{2a}\beta_{20} \\
   \beta_{02}=\xi^{a}\beta_{10}=\xi^{2a}\beta_{21}\\
   \beta_{012}=\xi^a\beta_{012}.
\end{array}
\end{equation}

To compute the relations between the coefficients for the action of $\iota$ we apply it to the equations of $\Phi_2|_{\UU_0}$, and get
\begin{equation}\label{conditionsiota}
    \begin{array}{rcl}
   \beta_1 &=& (-1)^b\beta_2 \\
   \beta_{01}&=&(-1)^b\beta_{02} \\
   \beta_{12}&=&(-1)^b\beta_{21} \\
    \beta_{20}&=&(-1)^b\beta_{02} \\
    \beta_{012}& = &(-1)^b\beta_{012}.
\end{array}
\end{equation}

For the action of $\tau\in\mu_3$, using Equations (\ref{equations}), we have the relations
$$c_0=\xi^{-c+2}\tau(c_0),\,\, c_1=\xi^{-c}\tau(c_1),\,\,
c_2=\xi^{-c+1}\tau(c_2),\,\, c_3=\xi^{-c+2}\tau(c_3),$$
from where we conclude that there are three options.
\begin{equation}\label{conditionstau}
    \begin{array}{ll}
        \beta_{01}=\beta_{02}=\beta_{10}=\beta_{12}=\beta_{20}=\beta_{21}=0,  & \text{ if $c=0$,} \\
        \beta_0=\beta_1=\beta_2=\beta_{012}=\beta_{10}=\beta_{21}=\beta_{02}=0, &  \text{ if $c=1$,} \\
        \beta_0=\beta_1=\beta_2=\beta_{012}=\beta_{01}=\beta_{12}=\beta_{20}=0, &  \text{ if $c=2$.} \\
    \end{array}
\end{equation}

Equating all the relations we conclude that if one of the integers $a,b,c$ is non zero (module the respective order), then the morphism equivariant for $H^e(1,3)$ is the zero map. 

\begin{proposition}\label{equivariantprop}
Given a covering homomorphism $\Phi_2\in\Hom\left(\Sym^2\Omega^1_{\PP^2},\Omega^1_{\PP^2}\right)$, coordinates $\{x_i\}$ for $\PP^2$ and a basis $\{y_i\}$ for $(\Omega^1_{\PP^2})^\vee$ such that $\Phi_2|_{\{x_2\neq 0\}}$ is of the form
\begin{equation*}
\begin{array}{rcl}
y_0^2 & = & c_0y_0+c_1y_1\\ 
y_0y_1 & = & -c_3y_3-c_0y_1 \\
y_1^2 & = & c_2y_2 + c_3y_1
\end{array}
\end{equation*}
with the $c_i$ as in Equations (\ref{equationsc}). Then $\Phi_2$ is equivariant for $H^e(1,3)$ acting as
\begin{equation*}\begin{array}{ccc}
\sigma(x_i,y_j) = (x_{i+1},y_{j+1}),\,\, &
\iota(x_i,y_j) = (x_{-i},y_{-j}),\,\, &
\tau(x_i,y_j) = (\xi^i x_i,\xi^{-i}y_i),
\end{array}\end{equation*}
if 
$\beta_0 = \beta_1 = \beta_2$ and $\beta_{01} = \beta_{02} = \beta_{10} = \beta_{12} = \beta_{20} = \beta_{21}=0.$
\end{proposition}

\begin{remark}
Taking the extension 
\begin{equation*}\begin{array}{cc}
\sigma(x_i,y_j) = (x_{i+1},y_{j+1}), & \,\,
\tau(x_i,y_j) = (\xi^i x_i,\xi^{c} \xi^{-i}y_i),  
\end{array}\end{equation*}
with $c\neq 0$ gives is an equivariant map for $H(1,3)$ but not for $H^e(1,3)$. These morphisms induce smooth isolated surfaces.
\end{remark}

\begin{remark}
 The choice of index in the coefficients $\beta_*$ was not random. If we assign to the coefficient $\beta_i$ the monomial $x_i^3$, to $\beta_{ij}$ the monomial $x_i^2x_j$ and to $\beta_{012}$ the monomial $x_0x_1x_2$, each of the actions $\sigma,\iota,\tau$ changes them exactly as if they are acting on those monomials.
 
 We have no explanation for this fact but it would be something worth further exploration, first studying if we have such a relation for $m>0$, and if we get some similar structure for covering homomorphisms in $\Hom\left(\Sym^2\Omega^1_{\PP^3},\Omega^1_{\PP^3}\right)$.
\end{remark}

%% file: CoordinateRing.tex
In this Section we describe the graded coordinate ring of $(S,\LL,\theta)$, where $S$ is an abelian \mbox{surface}, $\LL$ is an ample, symmetric line bundle that induces a polarisation of type $(1,3)$ and $\theta\colon G(\LL)\rightarrow H(1,3)$ is a theta-structure that induces an isomorphism between the representation $\widetilde\rho\colon G(\LL)\rightarrow \GL\left(H^0(S,\LL)^\vee\right)$ and the Schr\"odringer representation of $H(1,3)$ in $\GL(\operatorname{Map}(\ZZ_3,k^*))$.

Assuming that $(S,\LL)\ncong(E_1\times E_2,p_1^*\LL_1\otimes p_2^*\LL_2)$, where $E_i$ is an elliptic curve and $\LL_i$ a line bundle on $E_i$, then $|\LL|$ defines a covering map of degree $6$, $\varphi\colon S\rightarrow\PP^2$ such that
$$\varphi_*\OO_S=\OO_{\PP^2}\oplus \Omega^1_{\PP^2}\oplus \Omega^1_{\PP^2}\oplus \OO_{\PP^2}(-3).$$

The local structure of $\OO_S$ is determined by the following Theorem. We present a sketch of the proof taken from \cite[Thm. $3.2$, Thm. $3.3$]{CoverHomomorphisms}.

\begin{theorem}\label{abelianstructure}
 Let $(S,\LL)$ be a $(1,3)$-polarised abelian surface such that $(S,\LL)\ncong (E_1\times E_2,p_1^*\LL_1\otimes p_2^*\LL_2)$ and $\varphi\colon S\rightarrow \PP^2$ the covering map defined by the linear system $|\LL|$. Then
 \begin{enumerate}
     \item there is an embedding $S\hookrightarrow\AA\left((\Omega^1_{\PP^2}\oplus \Omega^1_{\PP^2})^\vee\right)$;
     \item a commutative and associative structure on $\varphi_*\OO_S$ determines and is determined by a cover homomorphism $\Phi\in\Hom\left((\Sym^2\Omega^1_{\PP^2})^{\oplus 3},(\Omega^1_{\PP^2})^{\oplus 2}\right)$;
     \item if $(y_0,y_1;z_0,z_1)$ is a local basis for $\Omega^1_{\PP^2}\oplus \Omega^1_{\PP^2}$ over the open set $\UU_2=\{x_2\neq 0\}$, then the local basis for $\Sym^2\left(\Omega^1_{\PP^2}\oplus \Omega^1_{\PP^2}\right)$ is 
     $$q=\left(
\begin{array}{ccccccccc}
 y_0^2 & y_0y_1 & y_1^2 & y_0z_1  & \frac{1}{2}(y_0z_1+y_1z_0) & y_1z_1 & z_0^2 & z_0z_1 & z_1^2
\end{array}\right),$$
    and $\Phi(q)=(\begin{array}{cccc}y_0 & y_1 & z_0 & z_1\end{array})C$, where $C$ is a $4\times 9$ matrix that decomposes as
    \begin{equation}\label{tracefreeglobal}
    C^t=\left(\begin{array}{cc}C_1  & C_0 \\ -C_2 & -C_1 \\ C_3   &  C_2
    \end{array}\right),\,\,
    C_i=\left(\begin{array}{cc}c_{i1}  & c_{i0} \\ -c_{i2}  & -c_{i1} \\ c_{i3}  & c_{i2}
    \end{array}\right),
    \end{equation}
    where the $c_{ij}\in\OO_{\PP^2}(\UU_2)$ satisfy the relations given by the spinor embedding of the orthogonal Grassmann variety $\oGr(5,10)$ in $\PP^{15}$.
 \end{enumerate}
\end{theorem}

\begin{proof}[Sketch of proof]
Using the adjunction formula and the fact that there exist a line bundle, $\MM=\OO_{\PP^2}(1)$, such that $\omega_S\cong \varphi^*\MM^{k_S}$ and $\omega_{\PP^2}\cong\MM^{k_{\PP^2}}$, for integers $k_S, k_{\PP^2}$, we conclude that $\varphi_*\OO_S$ and $\varphi_*\omega_S=\HHom(\varphi_*\OO_S,\omega_{\PP^2})$ are isomorphic as $\OO_S$-modules.

Picking one such isomorphism $\Psi\colon\varphi_*\OO_S\rightarrow\varphi_*\omega_S$, let $\eta:=\Psi\left(\varphi^\#(1_{\OO_{\PP^2}})\right)$ be the generator of $\varphi_*\omega_S$ as an $\OO_S$-module, and $e\in\OO_S$ the element such that $e\cdot\eta=\tfrac{1}{d}\operatorname{tr}$. Moreover, define an $\OO_{\PP^2}$-bilinear form $\langle,\rangle\colon\varphi_*\OO_S\times\OO_S\rightarrow\OO_{\PP^2}$ as 
$$ \langle a,b\rangle=\Psi(a)(b)=a\Psi\left(\varphi^\#(1_{\OO_{\PP^2}})\right)(b)=a\eta(b)=\eta(ab).$$

Then $\Omega^1_{\PP^2}\oplus\Omega^1_{\PP^2}$ is the orthogonal complement of $1_{\OO_S}$ and $e$. Furthermore, one can prove that when restricted to $\Omega^1_{\PP^2}\oplus\Omega^1_{\PP^2}$ the form $\langle,\rangle$ is nonsingular and $\eta|_{\Sym^2\left(\Omega^1_{\PP^2}\oplus\Omega^1_{\PP^2}\right)}$ is surjective. As $\OO_{\PP^2}\oplus\Omega^1_{\PP^2}\oplus\Omega^1_{\PP^2}$
are in the kernel of $\eta$, one concludes that $$\varphi_*\OO_S\subset\OO_{\PP^2}+\Omega^1_{\PP^2}\oplus\Omega^1_{\PP^2}+\Sym^2\left(\Omega^1_{\PP^2}\oplus\Omega^1_{\PP^2}\right).$$ Using the description of $\AA\left((\Omega^1_{\PP^2}\oplus\Omega^1_{\PP^2})^\vee\right)$ from Section \ref{Preliminaries} proves $(1)$. 

As $\Omega^1_{\PP^2}$ is a simple $\OO_{\PP^2}$-module
$$ 
\begin{array}{rcl}
    \Hom\left(\Omega^1_{\PP^2},\Omega^1_{\PP^2}\right)= k & \Leftrightarrow & \Hom\left(\Omega^1_{\PP^2}\otimes (\Omega^1_{\PP^2})^\vee,\OO_{\PP^2}\right)=k \\
    & \Leftrightarrow & \Hom\left(\Omega^1_{\PP^2}\otimes (\Omega^1_{\PP^2})^\vee\otimes \left(\Omega^1_{\PP^2}\bigwedge \Omega^1_{\PP^2}\right),\left(\Omega^1_{\PP^2}\bigwedge \Omega^1_{\PP^2}\right)\right)=k \\
    & \Leftrightarrow & \Hom\left(\Omega^1_{\PP^2}\otimes \Omega^1_{\PP^2},\OO_{\PP^2}(-3)\right)=k.
\end{array}
$$
From the decomposition $\Sym^2(\Omega^1_{\PP^2}\oplus\Omega^1_{\PP^2})=\left(\Sym^2 \Omega^1_{\PP^2}\right)^{\oplus 3}\oplus \OO_{\PP^2}(-3)$ we see that the ideal $\II_S\subset \OO_{\VV}$ is defined by a morphism
$$\Phi\in \Hom\left(\left(\Sym^2 \Omega^1_{\PP^2}\right)^{\oplus 3},\OO_{\PP^2}\oplus \Omega^1_{\PP^2}\oplus \Omega^1_{\PP^2}\right).$$
Then $(3)$ is a direct computation (see \cite[Thm. $3.3$]{CoverHomomorphisms}) and implies $(2)$ as decomposing $\Phi=\Phi_1\oplus\Phi_2$, where $\Phi_1$ is the morphism to $\Omega^1_{\PP^2}\oplus\Omega^1_{\PP^2}$ and $\Phi_2$ to $\OO_{\PP^2}$, then $\Phi_2$ is determined by $\Phi_1$.
\end{proof}

Using the previous Theorem and some computer algebra, one could use the local structure of a general covering homomorphism in $\Hom(\Sym^2\Omega^1_{\PP^2},\Omega^1_{\PP^2})$, see (\ref{equationsc}), to describe the coordinate ring $\RR(S,\LL)$. It would involve $15$ polynomials and $40$ variables.
By requiring the morphism $\varphi\colon S\rightarrow\PP^2$ to be equivariant for $H^e(1,3)$, as we will see in the next Theorem, the number of coefficients will become quite small and the moduli space easy to deal with.

\begin{theorem}\label{maintheorem}
 Let $S$ be an abelian surface over an algebraically closed field $k$ with characteristic different from $2$ and $3$, let $\LL$ an ample symmetric line bundle inducing a \mbox{polarisation} of type $(1,3)$ and let $\theta$ a canonical level structure. Then, for an open and dense subset of triples $(S,\LL,\theta)$, 
 the coordinate ring of $(S,\LL,\theta)$ can locally be written as
 $$\RR(S,\LL,\theta)|_{x_2\neq 0}=k[x_0,x_1,x_2^{\pm 1},y_0,y_1,z_0,z_1]/I,
 $$
 with the ideal $I$ generated by the entries of the vector 
 $$q-(\begin{array}{cccc}y_0 & y_1 & z_0 & z_1\end{array})C -D,$$
 where 
 \begin{itemize}
     \item $q$ and $C$ are as in Theorem $(\ref{abelianstructure})$;
     \item each block $C_i$ in $C$ has the entries 
      $$\begin{array}{rclrcl}
    c_{i0} & = & \beta_ix_2^2-\beta_i\frac{x_1^3}{x_2}, & 
    c_{i2} & = & -\frac{2}{3}\alpha_ix_0x_2-\beta_i\frac{x_0^2x_1}{x_2}, \\
    c_{i1} & = & -\frac{2}{3}\alpha_ix_1x_2-\beta_i\frac{x_0x_1^2}{x_2}, &
    c_{i3} & = &  \beta_ix_2^2-\beta_i\frac{x_0^3}{x_2};
 \end{array} 
 $$
    \item $D$ is a vector with entries
 $$D^t=\left(\begin{array}{l}
-2c_{11}^2 + 2c_{10}c_{12} + 2c_{01}c_{21} - c_{02}c_{20} - c_{00}c_{22} \\ 
-c_{10}c_{13} + c_{11}c_{12} - 2c_{02}c_{21} + c_{03}c_{20} + c_{01}c_{22} \\ 
2c_{11}c_{13} - 2c_{12}^2 - c_{03}c_{21} - c_{01}c_{23} + 2c_{02}c_{22} \\ 
-c_{01}c_{31} + c_{00}c_{32} + c_{11}c_{21} + c_{12}c_{20} - 2c_{10}c_{22} \\ 
\tfrac{1}{2}( -c_{00}c_{33} + c_{01}c_{32} - 5c_{12}c_{21} + c_{13}c_{20} + 4c_{11}c_{22}) \\ 
c_{01}c_{33} - c_{02}c_{32} + c_{13}c_{21} - 2c_{11}c_{23} + c_{12}c_{22} \\ 
2c_{11}c_{31} - c_{12}c_{30} - c_{10}c_{32} - 2c_{21}^2 + 2c_{20}c_{22} \\ 
c_{12}c_{31} + c_{10}c_{33} - 2c_{11}c_{32} - c_{20}c_{23} + c_{21}c_{22} \\ 
-c_{13}c_{31} - c_{11}c_{33} + 2c_{12}c_{32} + 2c_{21}c_{23} - 2c_{22}^2 
\end{array}\right);
$$
     \item and the coefficients $\alpha_i,\beta_i$ satisfy the equation
 \begin{equation}\label{moduliequation}
    3\alpha_0\beta_3-\alpha_1\beta_2+\alpha_2\beta_1-3\alpha_3\beta_0=0.
 \end{equation}
 \end{itemize}
\end{theorem}

\begin{proof}
Let $\sigma,\iota,\tau$ be lifts of $\sigma',\iota',\tau'$ that act on the variables $\{x_i,y_i,z_i\}_{0\leq i\leq 2}$ as
$$\begin{array}{rcl}
     \sigma(x_i,y_j,z_k) & = & (x_{i+1},y_{j+1},z_{k+1})\\
     \iota(x_i,y_j,z_k) & = & (x_{-i},y_{-j},z_{-k})\\
     \tau(x_i,y_j,z_k) & = & (\xi^i x_i,\xi^{-j} y_j,\xi^{-k} z_k).  \\
\end{array}
$$
Our goal is to prove that, independently of the choice of basis for $\left(\Omega_{\PP^2}^1\oplus\Omega^1_{\PP^2}\right)^\vee$, the action of $H^e(1,3)$ on $\RR(S,\LL,\theta)$ is generated by $\sigma, \iota, \tau$.

Let $\widetilde\sigma$ be a lift of $\sigma'$ that fixes the space generated by the polynomials $\sum x_iy_i, \sum x_iz_i$. Then 
$$\widetilde\sigma(x_i,y_j,z_k)=(x_{i+1},ay_{j+1}+bz_{j+1},cy_{k+1}+dz_{k+1}),$$
where $g=\left(\begin{array}{cc}
    a & b \\
    c & d
\end{array}\right)$ is a matrix in $\GL_2(k)$ such that $g^3=I_{2\times 2}$. As $k$ is an algebraically closed field with characteristic different from $3$, $g$ is diagonalisable. With a change of basis, we can consider $g=\operatorname{diag}(\xi^{m_\sigma},\xi^{n_\sigma})$, for $m_\sigma,n_\sigma\in \{0,1,2\}$.

Using such basis, the covering homomorphism $\Phi\in\Hom\left((\Sym^2\Omega_{\PP^2}^1)^{\oplus 3},(\Omega_{\PP^2}^1)^{\oplus 2}\right)$ in its local matrix form, we get that
\begin{equation}\label{conditionstildesigma}
\widetilde\sigma(\Phi)=\left(\begin{array}{cc}
     \xi^{-m_\sigma} \widetilde\sigma(C_1) & \xi^{m_\sigma+n_\sigma} \widetilde\sigma(C_0) \\
     -\xi^{-n_\sigma} \widetilde\sigma(C_2) & -\xi^{-m_\sigma} \widetilde\sigma(C_1) \\
     \xi^{m_\sigma+n_\sigma} \widetilde\sigma(C_3) & \xi^{-n_\sigma} \widetilde\sigma(C_2) \\
\end{array}\right).
\end{equation}
Then we have that the relations on the local coefficients defining each $C_i$ are given by $(\ref{conditionssigma})$ where $a\in\{m_\sigma,-m_\sigma-n_\sigma,n_\sigma\}$. 
Through computer algebra, see Appendix, we can prove that if $m_\sigma\neq n_\sigma$ then we are in one of the cases:
\begin{itemize}
    \item $C_0=\gamma C_3$, for a constant $\gamma\in k$, and $C_1=C_2=0$;
    \item $C_0 = C_2 = C_3 = 0$;
    \item $C_0 = C_1 = C_3 = 0$.
\end{itemize}    

The last two are degenerate. The first case has the property that the relations between the entries $\beta^i_*$ are the same for $C_0$ and $C_3$, which is a property valid for the case $m_\sigma=n_\sigma$. To be more specific, if $m_\sigma=n_\sigma$, from Section \ref{groupactions}, we have that for all $C_i$ their respective coefficients satisfy the following equalities
$$\begin{array}{c}
   \beta^i_0=\xi^{m}\beta^i_1=\xi^{-m}\beta^i_2 \\
   \beta^i_{01}=\xi^{m}\beta^i_{12}=\xi^{-m}\beta^i_{20} \\
   \beta^i_{10}=\xi^{m}\beta^i_{21}=\xi^{-m}\beta^i_{02} \\
   \beta^i_{012}=\xi^m\beta^i_{012}.
\end{array}$$

Then, using any basis, the local matrix of $\Phi$ is of the form
$$
\widetilde\sigma(\Phi)=\left(\begin{array}{cc}
     \xi^{-m_\sigma} \widetilde\sigma(\widetilde C_1) & \xi^{-m_\sigma} \widetilde\sigma(\widetilde C_0) \\
     -\xi^{-m_\sigma} \widetilde\sigma(\widetilde C_2) & -\xi^{-m_\sigma} \widetilde\sigma(\widetilde C_1) \\
     \xi^{-m_\sigma} \widetilde\sigma(\widetilde C_3) & \xi^{-m_\sigma} \widetilde\sigma(\widetilde C_2) \\
\end{array}\right),
$$
where each $\widetilde C_i$ is a linear combination of the $C_i$. The linear transformation $(y_i,z_i)\mapsto h(y_i,z_i)$, $h\in\GL_2(k)$, changes the matrix into
$$
\left(\begin{array}{ccc}
    h_{11}^2 & 2h_{11}h_{12} & h_{12}^2 \\
    h_{11}h{21} & h_{11}h_{22}+h_{12}h_{21} & h_{12}h_{22} \\
    h_{21}^2 & 2h_{21}h_{22} & h_{22}^2 \\
\end{array}\right)^{-1}
\left(\begin{array}{cc}
      C_1 &  C_0 \\
     - C_2 & - C_1 \\
     C_3 & C_2 \\
\end{array}\right)
\left(\begin{array}{cc}
    h_{11} & h_{12} \\
    h_{21} & h_{22} 
\end{array}\right).
$$

In particular, the relations between the coefficients of the entries of each matrix are kept. This implies that we can consider $\sigma(x_i,y_j,z_k)=(x_{i+1},\xi^{m_\sigma} y_{j+1},\xi^{m_\sigma} z_{k+1})$, $m_\sigma\in \{0,1,2\}$, independently of the basis.

As $k=\overline{k}$ and the characteristic of $k$ is not $2$, we can pick a basis such that
$\iota(x_i,y_j,z_k)=(x_{-i},(-1)^{m_\iota} y_{-j},(-1)^{n_\iota}z_{-k})$. From the relations imposed by $\iota$, Equations (\ref{conditionsiota}), independently of the signs $m_\iota, n_\iota$, we conclude directly that $m_\sigma=0$.

Finally, picking a a basis such that the lift of $\tau'$ acts as 
$\tau(x_i,y_j,z_k)=(\xi^i x_{i},\xi^{m_\tau} \xi^{-j} y_{j},\xi^{n_\tau}\xi^{-k} z_{k})$, via Magma computations, that if $m_\tau\neq n_\tau$, then we are in one the three cases of $m_\sigma\neq n_\sigma$. Again this allows us to consider $m_\tau= n_\tau$. 

From the conditions between the $\beta^i_*$ in Equations (\ref{conditionstau}), if $m_\tau\neq 0$, then only the coefficients $\beta^i_{jk}$, with $j^2+k$ fixed modulo $3$, could be different from zero. But then, by the action of $\iota$, all the $\beta^i_*$ vanish.
We conclude that $m_\tau=n_\tau=0$. The signs of $\iota$ are both positive, each negative one implies that the morphism is zero and we end with a degenerate case with all the $C_i$ linearly dependent.

The structure of each block $C_i$ in the matrix $C$ comes from Proposition \ref{equivariantprop} and the relation between the coefficients is $$3\beta_0^0\beta^3_{012}-\beta^1_0\beta^2_{012}+\beta^2_0\beta^1_{012}-3\beta^3_0\beta^0_{012},$$
that is Equation $(\ref{moduliequation})$ in the variables $\alpha_i, \beta_j$.
Vector $D$ is taken from \cite[Thm. $3.3$]{CoverHomomorphisms}.
\end{proof}

We now have everything for our last result. Notice that the local equations have some similarity with the equations of a triple cover. As $\Omega^1_{\PP^2}$ is a simple module, every isomorphism of $R(S,\LL,\theta)$ that fixes $\PP^2$ is given by a a linear change of variables 
$$(y_i,z_i)\mapsto (y_i,z_i)g,\,\,\,\, 0\leq i\leq 2,$$
where $g\in\PGL_2(k)$. With this in mind, consider three points, $\{p_0, p_1, p_2\}$, in the affine plane $\AA^2=\Spec(k[y,z])$, with coordinates given by the equations
\begin{equation}
    \left\{\begin{array}{rcl}
        y^2 & = & \alpha_1y +\alpha_0z + 2(\alpha_1^2-\alpha_0\alpha_2) \\
        yz & = & -\alpha_2y-\alpha_1z -(\alpha_1\alpha_2-\alpha_0\alpha_3) \\
        z^2 & = & \alpha_3y +\alpha_2z +2(\alpha_2^2-\alpha_1\alpha_3).
    \end{array}\right.
\end{equation}
From \cite[Prop. $3.4$]{CoverHomomorphisms}, the barycenter of the points $\{q_i\}$ is the origin in $\AA^2$. This implies that the coordinates of two of the points determine the third one. 

Using a change of variables $(y,z)\mapsto (y,z)q$, with $g\in\GL_2(k)$, if the three points are distinct, we can change two of to any other two of our choice. By \cite[Lemma $4.5$]{triple}, the three points are distinct if 
\begin{equation}\label{branchlocusine}
    \alpha_0^2\alpha_3^2+4\alpha_0\alpha_2^3-3\alpha_1^2\alpha_2^2+4\alpha_1^3\alpha_3-6\alpha_0\alpha_1\alpha_2\alpha_3\neq 0.
\end{equation}

\begin{corollary}\label{maincorollary}
The moduli space of abelian surfaces with a polarisation of type $(1,3)$ and canonical level structure, $\mathcal{A}^{lev}_{(1,3)}$, is rational.
\end{corollary}

\begin{proof}
Let $\RR(S,\LL,\theta)_{(\alpha_i,\beta_i)}$ be a model of a $(1,3)$-abelian surface with a level structure of canonical type. Assuming that the $\alpha_i$ satisfy the inequality (\ref{branchlocusine}), using a linear change of variables
$$(y_i,z_i)\mapsto (y_i,z_i)g$$
we can set $(\beta_0,\beta_1,\beta_2,\beta_3)=(0,1,1,0)$. Equation (\ref{moduliequation}) becomes $\alpha_1=\alpha_2$.

Notice that with the change of variables we fixed three points in the affine plane. These points are fixed by the action of a representation of the group $S_3$ in $\GL_2(k)$ generated by the elements
$$s=\left(\begin{array}{cc}
    0 & 1 \\
    -1 & -1
\end{array}\right),\,\,\,
r=\left(\begin{array}{cc}
    0 & 1 \\
    1 & 0
\end{array}\right).$$
The action on the parameters $(\alpha_0,\alpha_1,\alpha_3)$ is the following
$$\begin{array}{rcl}
    s & \colon & (-\alpha_0+\alpha_3,-\alpha_0+\alpha_1,-\alpha_0) \\
    s^2 & \colon & (-\alpha_3,\alpha_1-\alpha_3,\alpha_0-\alpha_3) \\
    r & \colon & (\alpha_3,\alpha_1,\alpha_0) \\
    rs & \colon & (\alpha_0-\alpha_3,\alpha_1-\alpha_3,-\alpha_3) \\
    rs^2 & \colon & (-\alpha_0,-\alpha_0+\alpha_1,-\alpha_0+\alpha_3). \\
\end{array}$$
With a change of basis $\delta_0=\alpha_+\alpha_1+\alpha_3,\, \delta_1=-3\alpha_0+\alpha_1+\alpha_3,\, \delta_2=\alpha_0+\alpha_1-3\alpha_3$, the action of $S_3$ is just the standard action
$$s(\delta_i)=\delta_{i+1}, \,\,\, r(\delta_i)=\delta_{-i},$$
and we get that $\mathcal{A}^{lev}_{(1,3)}$ is birational to $k^3$ as $$k[\delta_0,\delta_1,\delta_2]^{S_3}\cong k[s_0,s_1,s_2],$$
where $s_i$ is the symmetric polynomial of degree $i$ in the variables $\delta_j$.
\end{proof}

\subsection{Connection with previous results}

Suppose that $S=E\times E$, where $E$ is an elliptic curve over $k$. Then the line bundle
$$\LL= \OO_S(E\times\{0\}+\{0\}\times E+A)$$ 
where $A$ is the antidiagonal in $E\times E$, defines a polarisation of type $(1,3)$ on $S$. Moreover, in this case the covering map $\varphi\colon S\rightarrow \PP^2$ is a Galois covering map and in \cite{BLange}, Birkenhake and Lange, describe the branch locus of this family. It is $3C\in \PP^2$, where $C$ is defined by the vanishing of the polynomial 
$$
\begin{array}{l}
(x_0^6+x_1^6+x_2^6)+2(2\lambda^3-1)(x_0^3x_1^3+x_1^3x_2^3+x_2^3x_0^3) \\
\hspace{15mm} -6\lambda^2(x_0^4x_1x_2+x_0x_1^4x_2+x_0x_1x_2^4)-3\lambda(\lambda^3-4)x_0^2x_1^2x_2^2.
\end{array}
$$

Using the description of $\RR(S,\LL,\theta)$ from Corollary \ref{maincorollary}, this one dimensional model is obtained when taking $\alpha_0=\alpha_3=0$, and $\lambda=\alpha_1$.

A two-dimensional family is described by Casnati in \cite{Casab}. Consider the covering maps $\varphi\colon S\rightarrow\PP^2$ that decompose as $S\xrightarrow{\rho} S'\xrightarrow{\zeta}\PP^2$, where $\rho$ and $\zeta$ are double and triple covers, respectively. Such surfaces are called \textit{bielliptic} as they have a nontrivial involution $j\colon S\rightarrow S$.

The surface $S'$ is a ruled surface with invariant $e(S')=-1$ over an elliptic curve $E$. As a triple cover over $\PP^2$ it satisfies $\zeta_*\OO_{S'}=\OO_{\PP^2}\oplus\Omega^1_{\PP^2}$. From this fact, $\RR(S,\LL,\theta)$ describes these surfaces when $\{\alpha_1=0\}$ is on the $S_3$-orbit of the parameters $(\alpha_0,\alpha_1,\alpha_3)$.

\subsection{Irregular surface}\label{irregular}

Assume that $(S,\LL)$ is a pair of a surface and an ample divisor for which $|\LL|$ is a base point free inducing a covering map $\varphi_{|\LL|}\colon S\rightarrow\PP^2$ such that 
$$\varphi_*\OO_S=\OO_{\PP^2}\oplus\Omega^1_{\PP^2}(-m_1)\oplus\Omega^1_{\PP^2}(-m_2)\oplus\OO_{\PP^2}(-m_1-m_2-3),$$
for a pair of non-negative integers $(m_1,m_2)$. Then one can prove a result in the line of Theorem \ref{abelianstructure}. In particular, a commutative and associative multiplication structure on $\varphi_*\OO_S$ is determined by a cover homomorphism in 
$$\Hom\left(\Sym^2\Omega^1(-2m_1,-m_1-m_2,-2m_2),\Omega^1_{\PP^2}(-m_1,-m_2)\right),$$
where for a sheaf $\FF$ and integers $\{m_i\}_{1\leq i\leq r}$, we are denoting $\bigoplus_i\FF(m_i)$ by $\FF(m_1,\dots,m_r)$.

Each summand in these homomorphisms is described in Section \ref{triplecoverhomo}, and the local structure is determined by the relations of the embedding of the $\oGr(5,10)$. Then it is a computational problem. We write the equations for a irregular surface polarised by the canonical line bundle.

\begin{theorem}
 Let $(S,K_S)$ be a smooth surface with invariants $p_g=3$, $q=1$, $K_S^2=6$. If the linear system $|K_S|$ induces a covering map $\varphi\colon S\rightarrow \PP^2$ satisfying
 $$\varphi_*\OO_S=\OO_{\PP^2}\oplus\Omega^1_{\PP^2}\oplus\Omega^1_{\PP^2}(-1)\oplus\OO_{\PP^2}(-4),$$
 then there is a component of the moduli space of such surfaces whose canonical ring is described by a covering map in
    $$\Hom\left(\Sym^2\Omega^1(0,-1,-2),\Omega^1_{\PP^2}(0,-1)\right).$$
 For this component the canonical ring can locally be written as
 $$\RR(S,K_S)|_{\UU_2}=k[x_0,x_1,x_2^{\pm 1},y_0,y_1,z_0,z_1]/I,$$
 with $\deg(x_i,y_j,z_k)=(1,2,3)$, and the ideal $I$ is generated by the entries of
  $$q-(\begin{array}{cccc}y_0 & y_1 & z_0 & z_1\end{array})C -D,$$
  where 
  \begin{itemize}
      \item $q=\left(\begin{array}{ccccccccc}
 y_0^2 & y_0y_1 & y_1^2 & y_0z_1  & \frac{1}{2}(y_0z_1+y_1z_0) & y_1z_1 & z_0^2 & z_0z_1 & z_1^2
\end{array}\right)$;
        \item \begin{equation*}\label{tracefreeglobal}
    C^t=\left(\begin{array}{cc}C_1  & 0 \\ 0 & -C_1 \\ C_3   &  0
    \end{array}\right),\,\,
    C_i=\left(\begin{array}{cc}c_{i1}  & c_{i0} \\ -c_{i2}  & -c_{i1} \\ c_{i3}  & c_{i2}
    \end{array}\right),
    \end{equation*} 
    with the $c_{ij}$ satisfying the relation
    $$c_{13}c_{30}-3c_{12}c_{31}+3c_{11}c_{32}-c_{10}c_{33}=0;$$
        \item $D$ as described in Theorem \ref{maintheorem}.
  \end{itemize}
\end{theorem}

The construction is in all identical to the construction of the abelian surface in Theorem \ref{abelianstructure}. 
Notice that $\Hom\left(\Sym^2\Omega^1(-2),\OO_{\PP^2}(-4)\right)\neq 0$. That is the reason why such model only describes a component of the moduli space of these surfaces as claimed. 

\begin{proof}
For degree reasons $C_0=0$. Then, computing the relations between the $c_{ij}$ with all $c_{0j}=0$ we get
$$
\left\{\begin{array}{l}
     \operatorname{rk} \left(\begin{array}{cccc}
         c_{10} & c_{11} & c_{12} & c_{13} \\
         c_{20} & c_{21} & c_{21} & c_{23}  
     \end{array}\right)=0\\
     c_{13}c_{30}-3c_{12}c_{31}+3c_{11}c_{32}-c_{10}c_{33}=0 \\
     c_{23}c_{30}-3c_{22}c_{31}+3c_{21}c_{32}-c_{20}c_{33}=0.
\end{array}\right.
$$

Assuming that the $c_{1j}$ have no common divisor, there is $l\in \OO_{\PP^2}(\UU_2)$ such that $c_{2j}=lc_{1j}$. Then, with a change of basis $z_i\mapsto z_i+ly_i$, we get $C_2=0$ and we are done.
\end{proof}

%% file: Appendix.tex
To simplify notation we denote the matrices $(C_0,C_1,C_2,C_3)$ by $(B,A,D,C)$ (the notation used by Miranda for the entries defining a triple cover, see \cite{triple}), the entries $c_{ij}$ by the corresponding $A_j, B_j, C_j, D_j$ and the coefficients $\beta^i_*$ by $a_*, b_*,c_*, d_*$, respectively.

To prove that $m_\sigma=n_\sigma$ just run the code below running $m,n$ over all pairs in $\{0,1,2\}\times\{0,1,2\}$ with $m\neq n$, and choosing one of the variables $w_0, w_{01}, w_{10}$, where $w\in \{a,b,c,d\}$, different from zero.

\begin{verbatim}
K:=Rationals();
R<x>:=PolynomialRing(K);
K<e>:=ext<K|x^2+x+1>;

RR<a0,a1,a2,a01,a12,a20,a10,a21,a02,a012,
b0,b1,b2,b01,b12,b20,b10,b21,b02,b012,
c0,c1,c2,c01,c12,c20,c10,c21,c02,c012,
d0,d1,d2,d01,d12,d20,d10,d21,d02,d012>:=PolynomialRing(K,40);
R<x0,x1,x2>:=PolynomialRing(RR,3);

n:=*;m:=*;
*0:=1; or *01:=1; or *10:=1; 

a1:=(e^(-m))*a0; a2:=(e^(-m))^2*a0;
a12:=(e^(-m))*a01; a20:=(e^(-m))^2*a01;
a21:=(e^(-m))*a10; a02:=(e^(-m))^2*a10;

b1:=(e^(m+n))*b0; b2:=(e^(m+n))^2*b0;
b12:=(e^(m+n))*b01; b20:=(e^(m+n))^2*b01;
b21:=(e^(m+n))*b10; b02:=(e^(m+n))^2*b10;

c1:=(e^(m+n))*c0; c2:=(e^(m+n))^2*c0;
c12:=(e^(m+n))*c01; c20:=(e^(m+n))^2*c01;
c21:=(e^(m+n))*c10; c02:=(e^(m+n))^2*c10;

d1:=(e^(-n))*d0; d2:=(e^(-n))^2*d0;
d12:=(e^(-n))*d01; d20:=(e^(-n))^2*d01;
d21:=(e^(-n))*d10; d02:=(e^(-n))^2*d10;

A0 := a1*x2^3 - a12*x1*x2^2 + a21*x1^2*x2  - a2*(x1^3); 
A1 := (1/3)*a10*x2^3 - (2/3)*a012*x1*x2^2 + (1/3)*a20*x1^2*x2
- (1/3)*a12*x0*x2^2 + (2/3)*a21*x0*x1*x2 - a2*x0*x1^2;
A2 := (1/3)*a01*x2^3 - (1/3)*a02*x1*x2^2 - (2/3)*a012*x0*x2^2
+ (2/3)*a20*x0*x1*x2 + (1/3)*a21*x0^2*x2 - a2*x0^2*x1;
A3 := a0*x2^3  -  a02*x0*x2^2 +  a20*x0^2*x2 - a2*x0^3;
 
B0 := b1*x2^3 - b12*x1*x2^2 + b21*x1^2*x2  - b2*(x1^3); 
B1 := (1/3)*b10*x2^3 - (2/3)*b012*x1*x2^2 + (1/3)*b20*x1^2*x2
- (1/3)*b12*x0*x2^2 + (2/3)*b21*x0*x1*x2 - b2*x0*x1^2;
B2 := (1/3)*b01*x2^3 - (1/3)*b02*x1*x2^2 - (2/3)*b012*x0*x2^2
+ (2/3)*b20*x0*x1*x2 + (1/3)*b21*x0^2*x2 - b2*x0^2*x1;
B3 := b0*x2^3  -  b02*x0*x2^2 +  b20*x0^2*x2 - b2*x0^3;
 
C0 := c1*x2^3 - c12*x1*x2^2 + c21*x1^2*x2  - c2*(x1^3); 
C1 := (1/3)*c10*x2^3 - (2/3)*c012*x1*x2^2 + (1/3)*c20*x1^2*x2
- (1/3)*c12*x0*x2^2 + (2/3)*c21*x0*x1*x2 - c2*x0*x1^2;
C2 := (1/3)*c01*x2^3 - (1/3)*c02*x1*x2^2 - (2/3)*c012*x0*x2^2
+ (2/3)*c20*x0*x1*x2 + (1/3)*c21*x0^2*x2 - c2*x0^2*x1;
C3 := c0*x2^3  -  c02*x0*x2^2 +  c20*x0^2*x2 - c2*x0^3;
 
D0 := d1*x2^3 - d12*x1*x2^2 + d21*x1^2*x2  - d2*(x1^3); 
D1 := (1/3)*d10*x2^3 - (2/3)*d012*x1*x2^2 + (1/3)*d20*x1^2*x2 
- (1/3)*d12*x0*x2^2 + (2/3)*d21*x0*x1*x2 - d2*x0*x1^2;
D2 := (1/3)*d01*x2^3 - (1/3)*d02*x1*x2^2 - (2/3)*d012*x0*x2^2 
+ (2/3)*d20*x0*x1*x2 + (1/3)*d21*x0^2*x2 - d2*x0^2*x1;
D3 := d0*x2^3  -  d02*x0*x2^2 +  d20*x0^2*x2 - d2*x0^3;
 
oGr:=[
3*B3*A0 - 3*B0*A3 - B2*A1 + B1*A2,
3*B3*D0 - 3*B0*D3 - B2*D1 + B1*D2, 
3*B2*C0 - 3*B0*C2 - A2*D0 + A0*D2, 
3*B1*C0 - 3*B0*C1 - A1*D0 + A0*D1, 
9*B3*C0 - 9*B0*C3 - A2*D1 + A1*D2, 
3*B3*C2 - 3*B2*C3 - A3*D2 + A2*D3, 
3*B3*C1 - 3*B1*C3 - A3*D1 + A1*D3, 
3*D3*C0 - 3*D0*C3 - D2*C1 + D1*C2, 
3*A3*C0 - 3*A0*C3 - A2*C1 + A1*C2,
A3*D0 - A0*D3 - B2*C1 + B1*C2];

q1:=[(1-e^(-m))*a012,(1-e^(m+n))*b012,(1-e^(m+n))*c012,(1-e^(-n))*d012];
q2:=&cat[Coefficients(q):q in oGr];
q:=q1 cat q2;
I:=Ideal(q);
a0 in I; 
M:= MinimalBasis(I);
M;
\end{verbatim}

\noindent For the proof concerning $m_\tau, n_\tau$ each of the pieces of code below should trade the lines
\begin{verbatim}
n:=*;m:=*;
*0:=1; or *01:=1; or *10:=1;    
\end{verbatim}
in the code above.

\begin{verbatim}
m:=0;n:=1;
a01:=0;a10:=0;
b0:=0;b01:=0;b012:=0;
c0:=0;c01:=0;c012:=0;
d0:=0;d10:=0;d012:=0;

m:=1;n=2;
a0:=0;a10:=0;a012:=0;
b01:=0;b10:=0;
c01:=0;c10:=0;
d0:=0;d01:=0;d012:=0;

m:=1;n:=1;
a0:=0;a10:=0;a012:=0;
b0:=0;b10:=0;b012:=0;
c0:=0;c10:=0;c012:=0;
d0:=0;d10:=0;d012:=0;

m:=0;n:=0;
a01:=0;a10:=0;
b01:=0;b10:=0;
c01:=0;c10:=0;
d01:=0;d10:=0;
\end{verbatim}

%% file: AbSurf_Complete.bbl
\providecommand{\bysame}{\leavevmode\hbox to3em{\hrulefill}\thinspace}
\providecommand{\MR}{\relax\ifhmode\unskip\space\fi MR }
\providecommand{\MRhref}[2]{%
  \href{http://www.ams.org/mathscinet-getitem?mr=#1}{#2}
}
\providecommand{\href}[2]{#2}
\begin{thebibliography}{BLvS89}

\bibitem[Bei78]{Beilinson1}
A.~A. Beilinson, \emph{Coherent sheaves on {${\bf P}^{n}$} and problems in
  linear algebra}, Funktsional. Anal. i Prilozhen. \textbf{12} (1978), no.~3,
  68--69. \MR{509388}

\bibitem[BL94]{BLange}
Ch. Birkenhake and H.~Lange, \emph{{A family of abelian surfaces and curves of
  genus four}}, Manuscripta Math. \textbf{85} (1994), no.~3-4, 393--407.
  \MR{1305750 (95k:14064)}

\bibitem[BL04]{BLangebook}
Christina Birkenhake and Herbert Lange, \emph{{Complex abelian varieties}},
  second ed., {Grundlehren der Mathematischen Wissenschaften [Fundamental
  Principles of Mathematical Sciences]}, vol. 302, Springer-Verlag, Berlin,
  2004. \MR{2062673 (2005c:14001)}

\bibitem[BLvS89]{BLange4}
C.~Birkenhake, H.~Lange, and D.~van Straten, \emph{Abelian surfaces of type
  {$(1,4)$}}, Math. Ann. \textbf{285} (1989), no.~4, 625--646. \MR{1027763}

\bibitem[Cas99]{Casab}
Gianfranco Casnati, \emph{{The cover associated to a {$(1,3)$}-polarized
  bielliptic abelian surface and its branch locus}}, Proc. Edinburgh Math. Soc.
  (2) \textbf{42} (1999), no.~2, 375--392. \MR{1697405 (2000e:14078)}

\bibitem[Dia17]{CoverHomomorphisms}
Eduardo Dias, \emph{Construction of algebraic covers}, preprint
  arXiv:1709.03341 (2017).

\bibitem[GP98]{GrossPopescu}
Mark Gross and Sorin Popescu, \emph{Equations of {$(1,d)$}-polarized abelian
  surfaces}, Math. Ann. \textbf{310} (1998), no.~2, 333--377. \MR{1602020}

\bibitem[HM73]{Mumford5}
G.~Horrocks and D.~Mumford, \emph{A rank {$2$} vector bundle on {${\bf P}^{4}$}
  with {$15,000$}\ symmetries}, Topology \textbf{12} (1973), 63--81.
  \MR{0382279}

\bibitem[Man88]{15embedding}
Nicolae Manolache, \emph{Syzygies of abelian surfaces embedded in {${\bf
  P}^4({\bf C})$}}, J. Reine Angew. Math. \textbf{384} (1988), 180--191.
  \MR{929982}

\bibitem[Mir85]{triple}
Rick Miranda, \emph{{Triple covers in algebraic geometry}}, Amer. J. Math.
  \textbf{107} (1985), no.~5, 1123--1158. \MR{805807 (86k:14008)}

\bibitem[MS01]{MSchreyer7}
Nicolae Manolache and Frank-Olaf Schreyer, \emph{Moduli of {$(1,7)$}-polarized
  abelian surfaces via syzygies}, Math. Nachr. \textbf{226} (2001), 177--203.
  \MR{1839408}

\bibitem[Mum66]{Mumford1}
D.~Mumford, \emph{On the equations defining abelian varieties. {I}}, Invent.
  Math. \textbf{1} (1966), 287--354. \MR{0204427}

\bibitem[Pap04]{unproj1}
Stavros~Argyrios Papadakis, \emph{{Kustin-{M}iller unprojection with
  complexes}}, J. Algebraic Geom. \textbf{13} (2004), no.~2, 249--268.
  \MR{2047698 (2005d:13025)}

\bibitem[Pap06a]{PapadakisIII}
\bysame, \emph{Remarks on type {III} unprojection}, Comm. Algebra \textbf{34}
  (2006), no.~1, 313--321. \MR{2194769 (2006h:13031)}

\bibitem[Pap06b]{PapadakisII}
\bysame, \emph{Type {II} unprojection}, J. Algebraic Geom. \textbf{15} (2006),
  no.~3, 399--414. \MR{2219843 (2007c:14051)}

\bibitem[Pap07]{PapadakisGen}
\bysame, \emph{Towards a general theory of unprojection}, J. Math. Kyoto Univ.
  \textbf{47} (2007), no.~3, 579--598. \MR{2402516 (2009e:14019)}

\bibitem[PR04]{unproj2}
Stavros~Argyrios Papadakis and Miles Reid, \emph{{Kustin-{M}iller unprojection
  without complexes}}, J. Algebraic Geom. \textbf{13} (2004), no.~3, 563--577.
  \MR{2047681 (2005j:14068)}

\end{thebibliography}
